\documentclass{amsart}

\usepackage{amssymb,amscd}
\usepackage{epic, eepic, epsfig}
\usepackage{graphicx}
\newcounter{intro}

\newtheorem{thm}{Theorem}[section]
\newtheorem{lem}[thm]{Lemma}
\newtheorem{prop}[thm]{Proposition}
\newtheorem{cor}[thm]{Corollary}

\theoremstyle{remark}
\newtheorem{defi}[thm]{Definition}
\newtheorem{rem}[thm]{Remark}
\newtheorem{rems}[thm]{Remarks}

\numberwithin{equation}{section}   
\newcounter{counteroman}
\newenvironment{enumeroman}{\begin{list}{\roman{counteroman})}{\usecounter{counteroman}}}{\end{list}}

\newcommand{\cref}[1]{Corollary~\ref{#1}}

\newcommand{\R}{\mathbb{R}}
\newcommand{\C}{\mathbb{C}}

\newcommand{\N}{\mathbb{N}}
\newcommand{\Z}{\mathbb{Z}}
\newcommand{\bS}{\mathbb{S}}\newcommand{\bT}{\mathbb{T}}
\newcommand{\bD}{\mathbb{D}}\newcommand{\bB}{\mathbb{B}}
\newcommand{\cC}{\mathcal{C}}
\newcommand{\cH}{\mathcal{H}}\newcommand{\cA}{\mathcal{A}}
\newcommand{\cL}{\mathcal{L}}

\newcommand{\cD}{\mathcal{D}}

\let\priv=\smallsetminus

\DeclareMathOperator{\dive}{div}

\DeclareMathOperator{\capa}{cap}

\DeclareMathOperator{\Id}{Id}
\DeclareMathOperator{\ima}{Im}

\DeclareMathOperator{\inj}{inj}

\DeclareMathOperator{\ricci}{Ric}
\DeclareMathOperator{\Spec}{spec}
\DeclareMathOperator{\supp}{supp}

\DeclareMathOperator{\vol}{vol}
\DeclareMathOperator{\inte}{int}
\DeclareMathOperator{\grad}{\overrightarrow{grad}}

\begin{document}
\title{$L^2$ harmonics forms on non compact manifolds.}
\author{Gilles Carron}
\address{
Departement de Mathematiques\\
Universite de Nantes\\
2 rue de la Houssiniere, BP 92208\\
44322 Nantes Cedex 03, France\\ }
\email{Gilles.Carron\@math.univ-nantes.fr}
\maketitle
The source of these notes is a series of lectures given at the CIMPA's summer school 
"Recent Topics in Geometric Analysis". I want to thank the organizers of this summer school :
Ahmad El Soufi and Mehrdad Shahsahani and I also want to thank 
Mohsen Rahpeyma who solved many delicate problems.

Theses notes aimed to give an insight into some links between $L^2$ cohomology, $L^2$ harmonics forms, the
topology and the geometry of complete Riemannian manifolds. This is not a survey but a choice of few topics
in a very large subject.

The first part can be regard as an introduction ; we define the space of $L^2$ harmonics forms, of $L^2$
cohomology. We recall the theorems of Hodge and de Rham on compact Riemannian manifolds. However the
reader is assumed to be familiar with the basic of Riemannian geometry and with Hodge theory.

According to J. Roe (\cite{Roe}) and following the classification of Von Neumann algebra, we can classify problems
on $L^2$ harmonics forms in three types. The first one (type I) is the case where the space of harmonics
$L^2$ forms has finite dimension, this situation is the nearest to the case of compact manifolds. The second
 (type II) is the case where the space of harmonics
$L^2$ forms has infinite dimension but where we have a "renormalized" dimension for instance when a discrete
group acts cocompactly by isometry on the manifold; a good reference is the book of W. Lueck (\cite{Lu}) and the seminal
paper of M. Atiyah (\cite{at}). The third type (type III) is the case where no renormalization procedure
is available to define a kind of dimension of the space $L^2$ harmonics forms. Here we consider only the type I
problems and at the end of the first part, we will prove a result of J. Lott which says that the finiteness of
the dimension of the space of $L^2$ harmonics forms depends only on the geometry at infinity. 

Many aspects\footnote{almost all in fact !} of $L^2$ harmonics forms will not be treated here : for instance
we will not describe the important problem of the $L^2$ cohomology of locally symmetric spaces, and also we
will not speak on the pseudo differential approach developped by R. Melrose and his school.  However the reader
will find at the end of this first chapter a list of some interesting results on the topological
interpretation of the space of $L^2$ harmonics forms.

In the second chapter, we are interested in the space of harmonic $L^2$ $1-$forms. This space contains the 
differential of harmonic functions with $L^2$ gradient. We will not speak of the endpoint result of A.
Grigory'an (\cite{Gri, Gri0}) but we have include a study of P.Li and L-F. Tam (\cite{LT}) and of A. Ancona
(\cite{Ancona}) on
 non parabolic ends. In this chapter, we will also study the case of Riemannian surfaces where this space
 depends only on the complex structure.
 
 The last chapter focuses on the $L^2$ cohomology of conformally compact manifolds. The result is due to R.
 Mazzeo (\cite{M}) and the proof present here is the one of N. Yeganefar (\cite{Y1}) who used an
 integration by parts formula due to
 H.Donnelly and F.Xavier (\cite{DX}).

\tableofcontents

\section{A short introduction to $L^2$ cohomology}
In this first chapter, we introduce the main definitions and prove some
preliminary results.
\subsection{Hodge and de Rham 's theorems}
\subsubsection{de Rham 's theorem} Let $M^n$ be a smooth manifold of dimension $n$, we denote by 
$C^\infty(\Lambda^kT^*M)$ the space of smooth differential $k-$forms on $M$ and by
$C_0^\infty(\Lambda^kT^*M)$ the subspace of $C^\infty(\Lambda^kT^*M)$ formed by forms with compact support;
in local coordinates $(x_1,x_2,...,x_n)$, an element $\alpha \in C^\infty(\Lambda^kT^*M)$ has the following
expression
$$\alpha=\sum_{I=\{i_1<i_2<...<i_k\} }  \alpha_I dx_{i_1}\wedge dx_{i_2}\wedge ...\wedge dx_{i_k}
=\sum_{I }  \alpha_I dx_{I}$$
where $\alpha_I$ are smooth functions of $(x_1,x_2,...,x_n)$. The exterior differentiation
is a differential operator 
$$d\,:\, C^\infty(\Lambda^kT^*M)\rightarrow C^\infty(\Lambda^{k+1}T^*M)\, ,$$
locally we have
$$d\left(\sum_{I}  \alpha_{I} dx_{I}\right)=\sum_{I} d\alpha_{I}\wedge dx_{I}.$$

This operator satisfies $d\circ d=0$, hence the range of $d$ is included in the kernel of $d$.

\begin{defi}
The $k^{\rm th}$ de Rham's cohomology group of $M$ is defined by 
$$H^k_{dR}(M)=\frac{\left\{\alpha\in
C^\infty(\Lambda^kT^*M),d\alpha=0\right\}}{dC^\infty(\Lambda^{k-1}T^*M)}.$$
\end{defi}

These spaces are clearly diffeomorphism invariants of $M$, moreover the deep theorem of G. de Rham says that
these spaces are isomorphic to the real cohomology group of $M$, there are in fact
homotopy invariant of $M$ :
\begin{thm}
$$H^k_{dR}(M)\simeq H^k(M,\R).$$
\end{thm}
From now, we will suppress the subscript $dR$ for the de Rham's cohomology.
We can also define the de Rham's cohomology with compact support. 
\begin{defi}
The $k^{\rm th}$ de Rham's cohomology group with compact support of $M$ is defined by 
$$H^k_{0}(M)=\frac{\left\{\alpha\in
C_0^\infty(\Lambda^kT^*M),d\alpha=0\right\}}{dC_0^\infty(\Lambda^{k-1}T^*M)}.$$
\end{defi}

These spaces are also isomorphic to the real cohomology group of $M$ with compact support.
When $M$ is the interior of a compact manifold $\overline{M}$
 with compact boundary $\partial \overline{M}$
 $$M=\overline{M}\setminus\partial \overline{M},$$
 then $H^k_0(M)$ is isomorphic to the relative cohomology group of $\overline{M}$ :
 
 $$H_0^k(M)=H^k(\overline{M},\partial \overline{M}):=\frac{\left\{\alpha\in
C^\infty(\Lambda^kT^*\overline{M}),d\alpha=0, \iota^*\alpha=0\right\}}
{\left\{d\beta, \beta\in C^\infty(\Lambda^{k-1}T^*\overline{M}) \mbox{ and }\iota^*\beta=0\right\}}$$
where $\iota\,:\, \partial \overline{M}\rightarrow  \overline{M}$ is the inclusion map.

 \subsubsection{Poincar\'e duality}\label{PD} When we assume that $M$ is oriented \footnote{It is not a serious
 restriction we can used cohomology with coefficient in the orientation bundle.}
 the bilinear map 
 \begin{equation*}\begin{split}
 H^k(M) \times H_0^{n-k}(M)&\rightarrow \R\\
 ([\alpha],[\beta])&\mapsto \int_M \alpha\wedge \beta:=I([\alpha],[\beta])\\
 \end{split}\end{equation*}
 is well defined, that is to say $I([\alpha],[\beta])$ doesn't depend on the
 choice of representatives in the
 cohomology classes $[\alpha]$ or $[\beta]$ (this is an easy application of the Stokes formula). Moreover
this bilinear form provides an isomorphism between 
$H^k(M)$ and $\left(H_0^{n-k}(M)\right)^*$. In particular when $\alpha\in C^\infty(\Lambda^kT^*M)$ is closed
($d\alpha=0$) and satisfies that
$$\forall  [\beta]\in H_0^{n-k}(M),\ \int_M \alpha\wedge \beta=0$$ then there exists $\gamma\in
C^\infty(\Lambda^{k-1}T^*M)$ such that
$$\alpha=d\gamma.$$
 
\subsubsection{$L^2$ cohomology}\label{L2coho}

\noindent\ref{L2coho}{\it \,.a) The operator $d^*$.} We assume now that $M$ is endowed with a Riemannian metric $g$, we can
define the space $L^2(\Lambda^kT^*M)$ whose elements have locally the following
expression
$$\alpha=\sum_{I=\{i_1<i_2<...<i_k\}}  \alpha_{I} dx_{i_1}\wedge dx_{i_2}\wedge ...\wedge dx_{i_k}$$
where $\alpha_{I}\in L^2_{loc}$ and globally we have
$$\|\alpha\|_{L^2}^2:=\int_M |\alpha(x)|_{g(x)}^2 d\vol_g(x)<\infty.$$ 
The space $L^2(\Lambda^kT^*M)$ is a Hilbert space with scalar product :
$$\langle\alpha,\beta\rangle=\int_M (\alpha(x),\beta(x))_{g(x)} d\vol_g(x).$$
We define the formal adjoint of $d$ :
$$d^*\,:\,C^\infty(\Lambda^{k+1}T^*M)\rightarrow C^\infty(\Lambda^{k}T^*M)$$
by the formula 
\begin{equation*}\begin{split}
 &\forall \alpha \in
C_0^\infty(\Lambda^{k+1}T^*M)\mbox{ and } \beta\in C_0^\infty(\Lambda^{k}T^*M), \\
&\ \ \ \ \ \ \ \ \ \ \ \  \langle d^*\alpha,\beta\rangle=\langle\alpha,d\beta\rangle\ .\\
\end{split}\end{equation*}

When $\nabla$ is the Levi-Civita connexion of $g$, we can give local expressions for the operators $d$ and
$d^*$ : let $(E_1,E_2,..., E_n)$ be a local orthonormal frame and let 
$(\theta^1,\theta^2,...,\theta^n)$ be its dual frame :
$$\theta^i(X)=g(E_i,X)$$
then 
\begin{equation}\label{expd}
d\alpha=\sum_{i=1}^n\theta^i\wedge\nabla_{E_i}\alpha,
\end{equation}
and 
\begin{equation}\label{expd*}
d^*\alpha=-\sum_{i=1}^n\inte_{E_i}\left(\nabla_{E_i}\alpha\right),
\end{equation}
 where we have denote by 
 $\inte_{E_i}$ the interior product with the vector field $E_i$.
 \\
\noindent\ref{L2coho}{\it \,.b) $L^2$ harmonic forms.}
 We consider the space of $L^2$ closed forms :
 $$Z^k_2(M)=\{\alpha\in L^2(\Lambda^kT^*M),\ d\alpha=0\}$$
 where it is understood that the equation $d\alpha=0$ holds weakly that is to say
 $$\forall \beta\in C_0^\infty(\Lambda^{k+1}T^*M), \ \langle\alpha,d^*\beta\rangle =0.$$
That is we have :
 $$Z^k_2(M)=\left(d^*C^\infty(\Lambda^{k+1}T^*M)\right)^\perp\, ,$$
 hence $Z^k_2(M)$ is a closed subspace of $L^2(\Lambda^kT^*M)$.
 We can also define
 \begin{equation*}
 \begin{split}
 \cH^k(M)&=\left(d^*C_0^\infty(\Lambda^{k+1}T^*M)\right)^\perp\cap
 \left(dC_0^\infty(\Lambda^{k-1}T^*M)\right)^\perp\\
 &=Z^k_2(M)\cap\{\alpha\in L^2(\Lambda^kT^*M),\ d^*\alpha=0\}\\
 &=\{\alpha\in L^2(\Lambda^kT^*M),\ d\alpha=0\mbox{ and } d^*\alpha=0\}.\\
 \end{split}\end{equation*} 
 Because the operator $d+d^*$ is elliptic, we have by elliptic regularity :
 $\cH^k(M) \subset C^\infty(\Lambda^kT^*M)$.
 We also remark that by definition we have
 \begin{equation*}
 \begin{split}
\forall \alpha\in C_0^\infty(\Lambda^{k-1}T^*M),&\ \forall \beta\in C_0^\infty(\Lambda^{k+1}T^*M)\\
 \langle d\alpha,d^*\beta\rangle & =\langle dd\alpha, \beta\rangle=0\\
 \end{split}\end{equation*} 
 Hence $$dC_0^\infty(\Lambda^{k-1}T^*M)\perp d^*C_0^\infty(\Lambda^{k+1}T^*M)$$ and we get the Hodge-de
 Rham
 decomposition of $L^2(\Lambda^{k}T^*M)$
 \begin{equation}\label{decomposition}
 L^2(\Lambda^{k}T^*M)=\cH^k(M)\oplus \overline{dC_0^\infty(\Lambda^{k-1}T^*M)}\oplus
 \overline{ d^*C_0^\infty(\Lambda^{k+1}T^*M)}\ ,
 \end{equation} where the closures are taken for the $L^2$ topology.
 And also 
 \begin{equation}\label{harmo}
  \cH^k(M)\simeq \frac{Z^k_2(M)}{\overline{dC_0^\infty(\Lambda^{k-1}T^*M)}}.
  \end{equation}
 \noindent\ref{L2coho}{\it \,.c) $L^2$ cohomology:}
 We also define the (maximal) domain of $d$ by
 $$\cD^k(d)=\{\alpha\in L^2(\Lambda^kT^*M),\ d\alpha\in L^2\}$$
 that is to say 
$\alpha\in\cD^k(d)$ if and only if there is a constant $C$ such that 
$$\forall \beta\in C_0^\infty(\Lambda^{k+1}T^*M),\ |\langle\alpha,d^*\beta\rangle|\le C \|\beta\|_2.$$
In that case, the linear form $ \beta\in C_0^\infty(\Lambda^{k+1}T^*M)\mapsto\langle\alpha,d^*\beta\rangle$
extends continuously to $L^2(\Lambda^{k+1}T^*M)$ and there is $\gamma=:d\alpha$ such that
$$\forall\beta\in C_0^\infty(\Lambda^{k+1}T^*M),\  \langle\alpha,d^*\beta\rangle=\langle\gamma,\beta\rangle.$$
We remark that we always have $d\cD^{k-1}(d)\subset Z^k_2(M)$.
\begin{defi}
We define the $k^{\rm th}$ space of reduced $L^2$ cohomology by 
$$H^k_2(M)=\frac{ Z^k_2(M)}{\overline{d\cD^{k-1}(d)}}.$$
The $ k^{\rm th}$ space of non reduced $L^2$ cohomology is defined
by $${}^{nr}H^k_2(M)=\frac{ Z^k_2(M)}{d\cD^{k-1}(d)}.$$
 \end{defi}
 These two spaces coincide when the range of $d\,:\,\cD^{k-1}(d)\rightarrow L^2$ is closed; the first space
 is always a Hilbert space and the second is not necessary Hausdorff.
We also have $C_0^\infty(\Lambda^{k-1}T^*M)\subset \cD^{k-1}(d)$ hence we always get 
a surjective map :
$$\cH^k(M)\rightarrow H^k_2(M)\rightarrow\{0\}.$$
In particular any class of reduced $L^2$ cohomology contains a smooth representative.
\\
\noindent\ref{L2coho}{\it \,.d) Case of complete manifolds.}$\ $
The following result is due to Gaffney (\cite{gaff}, see also part 5 in \cite{wolf}) for a related result)
\begin{lem}\label{gaffney} Assume that $g$ is a complete Riemannian metric then
$$\overline{d\cD^{k-1}(d)}=\overline{dC_0^\infty(\Lambda^{k-1}T^*M)}.$$
\end{lem}
\proof We already know that $dC_0^\infty(\Lambda^{k-1}T^*M)\subset\cD^{k-1}(d)$, moreover using a
partition of unity and local convolution it is not hard to check that if $\alpha\in\cD^{k-1}(d)$ has compact
support then we can find a sequence
$(\alpha_l)_{l\in\N}$ of smooth forms with compact support such that
$$\|\alpha_l-\alpha\|_{L^2}+\|d\alpha_l-d\alpha\|_{L^2}\le 1/l.$$

So we must only prove that if $\alpha\in\cD^{k-1}(d)$ then we can build a sequence $ (\alpha_N)_N$ of elements of
$\cD^{k-1}(d)$ with compact support such that 
$$L^2-\lim_{N\to\infty}d\alpha_N=d\alpha .$$
We fix now an origin $o\in M$ and denote by $B(o,N)$ the closed geodesic ball of radius $N$ and centered at 
$o$, because $(M,g)$ is assumed to be complete we know that $B(o,N)$ is compact and
$$M= \cup_{N\in \N}B(o,N).$$
We consider $\rho\in C^\infty_0(\R_+)$ with $0\le \rho\le 1$ with support in $[0,1]$ such that 
$$\rho=1\mbox{ on } [0,1/2]$$
and we define
\begin{equation}\label{cutof}
\chi_N(x)=\rho\left(\frac{d(o,x)}{N}\right).
\end{equation}
Then $\chi_N$ is a Lipschitz function and is differentiable almost everywhere and
$$d\chi_N(x)=\rho'\left(\frac{d(o,x)}{N}\right) dr$$ where
$dr$ is the differential of the function $x\mapsto d(o,x)$.
Let $\alpha\in\cD^{k-1}(d)$ and define
$$\alpha_N=\chi_N\alpha,$$
the support of $\alpha_N$ is included in the ball of radius $N$ and centered at 
$o$ hence is compact. Moreover
we have 
$$\|\alpha_N-\alpha\|_{L^2}\le \|\alpha\|_{L^2(M\setminus B(o,N/2))}$$ hence
$$L^2-\lim_{N\to\infty}\alpha_N=\alpha$$
Moreover when $\varphi \in C_0^\infty(\Lambda^{k+1}T^*M)$ we have
\begin{equation*}\begin{split}
\langle \alpha_N,d^*\varphi\rangle&=\langle \alpha,\chi_Nd^*\varphi\rangle\\
&=\langle \alpha,d^*(\chi_N\varphi)\rangle+\langle \alpha,\inte_{\grad \chi_N}\varphi\rangle\\
&=\langle \chi_Nd\alpha+d\chi_N\wedge\alpha,\varphi\rangle\\
\end{split}\end{equation*}
Hence $\alpha_N \in\cD^{k-1}(d)$ and $$d\alpha_N=\chi_Nd\alpha+d\chi_N\wedge\alpha.$$
But for almost all $x\in M$, we have $|d\chi_N|(x)\le \|\rho'\|_{L^\infty}/ N$ hence
\begin{equation*}\begin{split}
\|d\alpha_N-d\alpha\|_{L^2}&
\le \frac{\|\rho'\|_{L^\infty}}{ N}\|\alpha\|_{L^2}+\|\chi_Nd\alpha-d\alpha\|_{L^2}\\
&\le\frac{\|\rho'\|_{L^\infty}}{ N}\|\alpha\|_{L^2}+\|d\alpha\|_{L^2(M\setminus B(o,N/2))}\ .\\
\end{split}\end{equation*}
Hence we have  build a sequence $\alpha_N$ of elements of
$\cD^{k-1}(d)$ with compact support such that 
$L^2-\lim_{N\to\infty}d\alpha_N=d\alpha.$
\endproof
A corollary of this lemma (\ref{gaffney}) and of (\ref{harmo}) is the following :
\begin{cor}When $(M,g)$ is a complete Riemannian manifold then the space of harmonic $L^2$ forms computes 
the reduced $L^2$ cohomology :
$$H^k_2(M)\simeq \cH^k(M).$$
\end{cor}
With a similar proof, we have another result :
\begin{prop}When $(M,g)$ is a complete Riemannian manifold then
$$\cH^k(M)=\{\alpha\in L^2(\Lambda^kT^*M), (dd^*+d^*d)\alpha=0\}.$$
\end{prop}
\proof Clearly we only need to check the inclusion :
$$\{\alpha\in L^2(\Lambda^kT^*M), (dd^*+d^*d)\alpha=0\}\subset\cH^k(M).$$
We consider again the sequence of cut-off functions $\chi_N$ defined previously in (\ref{cutof}).
Let $\alpha\in L^2(\Lambda^kT^*M)$ satisfying $(dd^*+d^*d)\alpha=0$ by elliptic regularity we know that
$\alpha$ is smooth. Moreover we have :
\begin{equation*}
\begin{split}
\left\|d(\chi_N\alpha)\right\|_{L^2}^2&=
\int_M\left[ |d\chi_N\wedge\alpha|^2+2\langle d\chi_N\wedge\alpha,\chi_N d\alpha\rangle+\chi_N^2|d\alpha|^2 \right]
d\vol_g\\
&=\int_M\left[ |d\chi_N\wedge\alpha|^2+\langle d\chi_N^2\wedge\alpha,d\alpha\rangle+\chi_N^2|d\alpha|^2 \right]
d\vol_g\\
&=\int_M\left[ |d\chi_N\wedge\alpha|^2+\langle d(\chi_N^2\alpha),d\alpha\rangle\right]
d\vol_g\\
&=\int_M\left[ |d\chi_N\wedge\alpha|^2+\langle \chi_N^2\alpha,d^*d\alpha\rangle\right]
d\vol_g\\
\end{split}
\end{equation*}
Similarly we get :
\begin{equation*}
\left\|d^*(\chi_N\alpha)\right\|_{L^2}^2=
\int_M\left[ |\inte_{\grad\chi_N}\alpha|^2+\langle \chi_N^2\alpha,dd^*\alpha\rangle\right]
d\vol_g
\end{equation*}
Summing these two equalities we obtain :
$$\left\|d^*(\chi_N\alpha)\right\|_{L^2}^2+\left\|d^*(\chi_N\alpha)\right\|_{L^2}^2=
\int_M |d\chi_N|^2|\alpha|^2d\vol_g\le \frac{\|\rho'\|^2_{L^\infty}}{N^2}\int_M |\alpha|^2d\vol_g.$$
Hence when $N$ tends to $\infty$ we obtain
$$\left\|d^*\alpha\right\|_{L^2}^2+\left\|d^*\alpha\right\|_{L^2}^2=0.$$
\endproof

This proposition has the consequence that on a complete Riemannian manifold harmonic $L^2$ functions are
closed hence locally constant. Another corollary is that the reduced $L^2$ cohomology of the Euclidean space
is trivial\footnote{This can also be proved with the Fourier transform.} :
\begin{cor}
$$H_2^k(\R^n)=\{0\}.$$
\end{cor}
\proof On the Euclidean space $\R^n$ a smooth $k$ form $\alpha$ can be expressed
as
$$\alpha =\sum_{I=\{i_1<i_2<...<i_k\}}  \alpha_{I} dx_{i_1}\wedge dx_{i_2}\wedge ...\wedge dx_{i_k}$$
and $\alpha$ will be a $L^2$ solution of the equation $$(dd^*+d^*d)\alpha=0$$ 
if and only if all the functions $\alpha_{I}$ are harmonic and $L^2$ hence zero because
the volume of $\R^n$ is infinite.\endproof

\begin{rem} When $(M,g)$ is not complete, we have not necessary equality between the space
$\cH^k(M)$ (whose elements are sometimes called harmonics fields) and the space of the
$L^2$ solutions of the equation $(dd^*+d^*d)\alpha=0$. For instance, on the interval $M=[0,1]$ the 
space $\cH^0(M)$ is the space of constant functions, whereas 
$L^2$ solutions of the equation $(dd^*+d^*d)\alpha=0$ are affine. More generally, on
a smooth compact connected manifold with smooth boundary endowed with a smooth Riemannian metric, then again 
$\cH^0(M)$ is the space of constant functions, whereas the space $\{f\in L^2(M), d^*df=0\}$ is the space of
harmonic $L^2$ function; this space is infinite dimensionnal when $\dim M>1$.
\end{rem}
\noindent\ref{L2coho}{\it e) Case of compact manifolds} The Hodge's theorem says that
 for compact manifold cohomology is computed with harmonic forms :
 \begin{thm} If $M$ is a compact Riemannian manifold without boundary
 then 
 $$H^k_ 2(M)\simeq \cH^k(M)\simeq H^k(M).$$
 \end{thm}
 
 When $M$ is the interior of a compact manifold $\overline{M}$
 with compact boundary $\partial \overline{M}$ and when $g$ extends to $\overline{M}$ (hence $g$ is
 incomplete) a theorem of P. Conner (\cite{conn}) states that
 $$H^k_ 2(M)\simeq H^k(M)\simeq\cH^k_{abs}(\overline{M})$$ where
  $$\cH^k_{abs}(M)=\{\alpha\in L^2(\Lambda^kT^*M),\ d\alpha=d^*\alpha=0\mbox{ and }
 \inte_{\vec\nu}\alpha=0 \mbox{ along }\partial \overline{M}\}$$
and $\vec \nu\,:\, \partial \overline{M}\rightarrow T\overline{M}$ is the inward unit normal vector field.
 In fact when $K\subset M$ is a compact subset of $M$ with smooth boundary and if $g$ is a complete Riemannian
 metric on $M$ then for $\Omega=M\setminus K$, we also have the equality
 $$H^k_ 2(\Omega)\simeq\cH^k_{abs}(\Omega)$$
 where if $\vec \nu\,:\, \partial \Omega\rightarrow TM$ is the inward unit normal vector field,
 we have also denoted 
 \begin{equation}\label{hodgeabs}\cH^k_{abs}(\Omega)=\{\alpha\in L^2(\Lambda^kT^*\Omega),\ d\alpha=d^*\alpha=0\mbox{ and }
 \inte_{\vec\nu}\alpha=0 \mbox{ along }\partial \Omega\}.\end{equation}
\subsection{Some general properties of reduced $L^2$ cohomology}
\subsubsection{a general link with de Rham's cohomology} 
We assume that $(M,g)$ is a complete Riemannian manifold, the following result is due to de Rham
(theorem 24 in \cite{dR})
\begin{lem} \label{dRham}Let $\alpha\in Z^k_2(M)\cap C^\infty(\Lambda^kT^*M)$ 
and suppose that $\alpha$ is zero in
$H^k_2(M)$ that is there is a sequence $\beta_l\in C^\infty_0(\Lambda^{k-1}T^*M)$ such that
$$\alpha=L^2-\lim_{l\to\infty} d\beta_l$$ then there is
$\beta\in C^\infty(\Lambda^{k-1}T^*M)$ such that 
$$\alpha=d\beta.$$
\end{lem}
In full generality, we know nothing about the behavior of $\beta$ at infinity.
\proof We can always assume that $M$ is oriented, hence by the Poincar\'e duality (\ref{PD}),
 we only need to show that
if $\psi\in C_0^\infty(\Lambda^{n-k}T^*M)$ is closed then
$$\int_M\alpha\wedge\psi =0\ .$$
But by assumption,
\begin{equation*}
\begin{split}
\int_M\alpha\wedge\psi &=\lim_{l\to \infty}\int_Md\beta_l\wedge\psi\\
&= \lim_{l\to \infty}\int_Md(\beta_l\wedge\psi) \mbox{ because }\psi\mbox{ is closed}\\
&=0. \\
\end{split}
\end{equation*}
Hence the result.\endproof

This lemma implies the following useful result which is due to M. Anderson (\cite{An}):
\begin{cor}\label{inject}
There is a natural injective map 
$$\ima\big( H^k_0(M)\rightarrow H^k(M)\big) \rightarrow H^k_2(M).$$
\end{cor}
\proof As a matter of fact we need to show that if $\alpha\in C_0^\infty(\Lambda^{k}T^*M)$ is closed  and zero
in the reduced $L^2$ cohomology then it is zero in usual cohomology: this is exactly the statement of the
previous lemma (\ref{dRham}).
\endproof
\subsubsection{Consequence for surfaces.}
These results have some implications for a complete Riemannian surface $(S,g)$ :
\begin{enumeroman}
\item If the genus of $S$ is infinite then the dimension of the space of $L^2$ harmonic $1-$forms is 
infinite.
\item If the space of $L^2$ harmonics $1-$forms is trivial then the genus of $S$ is zero and
$S$ is diffeomorphic to a open set of the sphere.\end{enumeroman}

As a matter of fact, a handle of $S$ is a embedding
$f\,:\, \bS^1 \times [-1,1]\rightarrow S$ such that if we denote $\cA=f(\bS^1 \times [-1,1])$ then
$S\setminus\cA$ is connected. 
\begin{center}\fbox{\includegraphics[height=6cm]{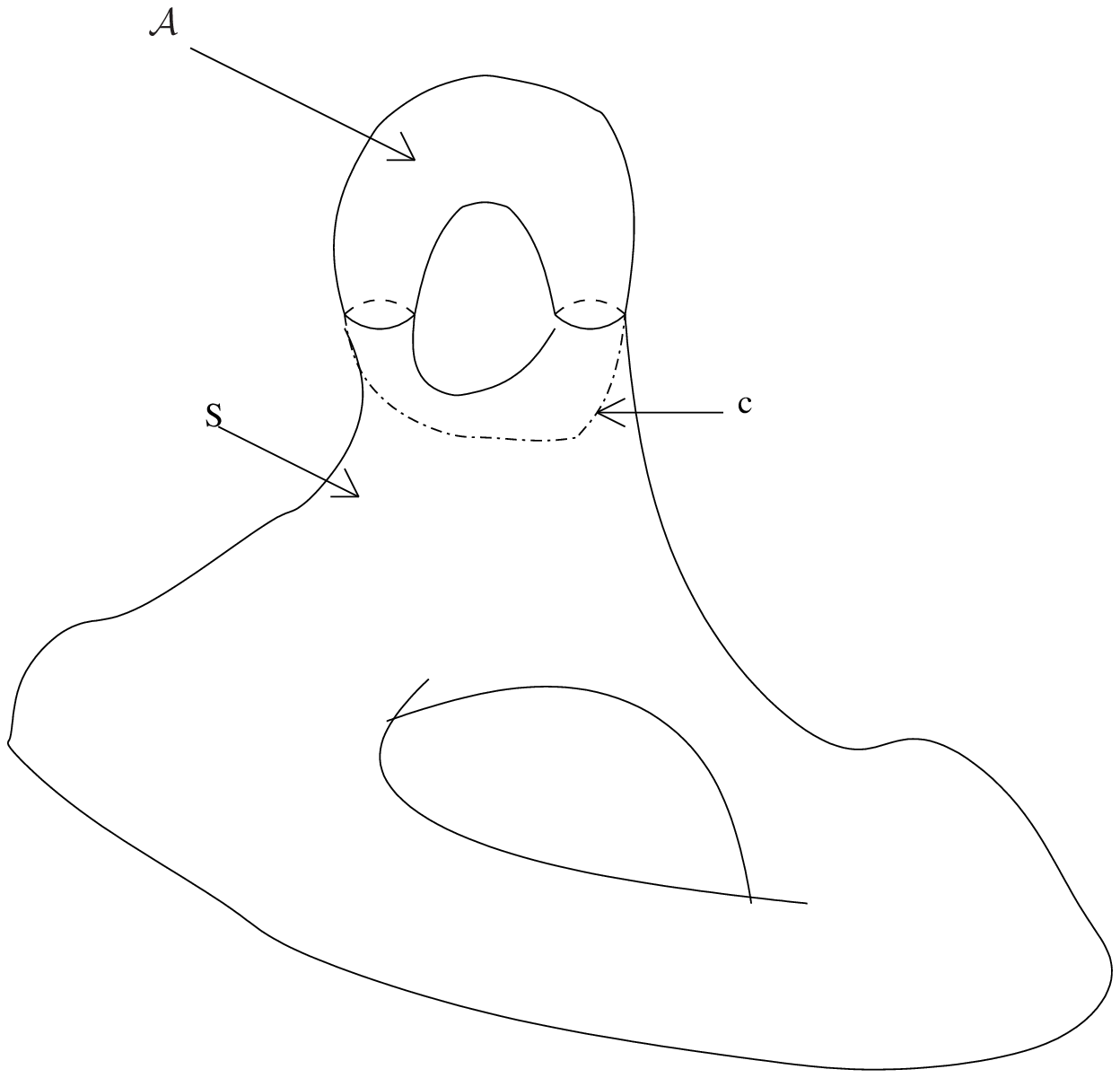}  }
\end{center}
We consider now a function $\rho$ on $\bS^1 \times [-1,1]$ depending only on the second variable such that
$$\rho(  \theta,t)=\left\{\begin{array}{ll}
1&\mbox{ when } t>1/2\\
0 &\mbox{ when } t<-1/2\\
\end{array}\right.$$
then $d\rho$ is a $1-$form with compact support in $\bS^1 \times]-1,1[$ and we can extend
$\alpha=\left(f^{-1}\right)^*d\rho$ to all $S$ ;  we obtain a closed $1$-form also denoted by $\alpha$ which
has  compact support in $\cA$ . Moreover
because $S\setminus\cA$ is connected, we can find a continuous path $c\,:\,[0,1]\rightarrow S\setminus\cA$ joining
$f(1,1)$ to $f(1,-1)$; we can defined the loop $\gamma$ given by
$$\gamma(t)=\left\{\begin{array}{ll}
f(1,t)&\mbox{ for } t\in [-1,1]\\
c(t-1) &\mbox{ for } t\in [1,2].\\
\end{array}\right.$$
It is easy to check that 
$$\int_\gamma \alpha=1\ ,$$ hence $\alpha$ is not zero in $H^1(S)$.
A little elaboration from this argument shows that 
$${\rm genus}(S)\le \dim\ima\big( H^1_0(M)\rightarrow H^1(M)\big).$$
\subsection{Lott's result} We will now prove the following result due to J.Lott (\cite{lottL2}):
\begin{thm} Assume that $(M_1,g_1)$ and $(M_2,g_2)$ are complete oriented manifold of dimension $n$ which are
isometric at infinity that is to say there are compact sets $K_1\subset M_1$ and $K_2\subset M_2$ such that
$(M_1\setminus K_1,g_1)$ and $(M_2\setminus K_2,g_2)$ are isometric. Then for $k\in[0,n]\cap\N$
$$\dim\cH^k(M_1,g_1)<\infty\Leftrightarrow \dim\cH^k(M_2,g_2)<\infty.$$
\end{thm}
We will give below the proof of this result, this proof contains many arguments which will be used and
refined in the next two lectures. In view of the Hodge-de Rham theorem and of J. Lott's result,
 we can ask the following very general questions :
\begin{enumerate}
\item What are the geometry at infinity iensuring the finiteness of the dimension of the spaces $H_2^k(M)$?

Within a class of Riemannian manifold having the same geometry at infinity :

\item What are the links of the spaces of reduced $L^2$ cohomology $H_2^k(M)$ 
with the topology of $M$ and with the geometry "at infinity" of $(M,g)$ ? 
\end{enumerate}

There is a lot of articles dealing with these questions, I mention some of them :
\begin{enumerate}
\item In the pioneering article of Atiyah-Patodi-Singer (\cite{APS}), the authors considered manifold
 with cylindrical end : that is to say
there is a compact $K$ of $M$ such that $M\setminus K$ is isometric to the Riemannian product 
$\partial K\times ]0,\infty[$. Then they show that the dimension of the
 space of $L^2$-harmonic forms is finite ; and that these spaces are isomorphic to the image of the relative
cohomology in the absolute cohomology. These results were used by Atiyah-Patodi-Singer in order to obtain a 
formula for the signature of compact Manifolds with boundary. 

\item In \cite{M, MP}, R. Mazzeo and R.Phillips give a cohomological interpretation of the space 
$\cH^k(M)$ for 
geometrically finite real hyperbolic manifolds. These manifolds can be compactified. 
They identify the reduced $L^2$ cohomology with the
cohomology of smooth differential forms satisfying certain boundary conditions.

\item The solution of the Zucker's conjecture by L.Saper-M.Stern  and E.Looijenga 
 (\cite{Loo},\cite{SS}) shows that the spaces of $L^2$ harmonic forms
on Hermitian locally symmetric space with finite volume are isomorphic to the
middle intersection cohomology of the Baily-Borel-Satake compactification of the manifold. An extension of this
result has been given by A.Nair and L.Saper (\cite{nair},\cite{saper1}). Moreover recently, L. Saper obtains the topological interpretation
of the reduced $L^2$ cohomology of any locally symmetric space with finite volume (\cite{saper2}). In that case the
finiteness of the dimension of the space of $L^2$ harmonics forms is due to A. Borel and H.Garland (\cite{BG}).
\item According to Vesentini (\cite{visentini})
 if $M$ is flat outside a compact set, the spaces $\cH^k(M)$ are finite dimensional. J. Dodziuk asked about
 the topological interpretation of the space $\cH^k(M)$ (\cite{Dod}). In this case, the answer has been given
 in \cite{Cargafa}.
 
 \item In a recent paper (\cite{HHM}) Tam\'as Hausel, Eugenie Hunsicker and Rafe Mazzeo obtain a topological
interpretation of the $L^2$ cohomology of complete Riemannian manifold whose geometry at infinity is fibred boundary 
and fibred cusp (see \cite{MM, Vaillant}).
 These results have important application concerning the Sen's conjecture (\cite{H, S}).
 
 \item In \cite{lottnegative}, J. Lott has shown that on a complete Riemannian manifold with finite volume and pinched negative
 curvature, the space of harmonic $L^2$ forms has finite dimension. N. Yeganefar obtains the topological
 interpretation of these spaces in two cases, first when the curvature is enough pinched (\cite{Y1}) and secondly
  when the metric is  K\"ahler (\cite{Y2}).
 \end{enumerate}
{\it Proof of J. Lott's result.} 
We consider $(M,g)$ a complete oriented Riemannian manifold and $K\subset M$ a compact subset with smooth
boundary and we let $\Omega=M \setminus K$ be the exterior of $K$, we are going to prove that 
$$\dim H^k_2(M)<\infty \Leftrightarrow \dim H^k_2(\Omega)<\infty\ ;$$
this result clearly implies Lott's result.

The co boundary map $b\,:\, H^k_2(\Omega)\rightarrow H^{k+1}(K,\partial K)$ is defined as follow:
let $c=[\alpha]\in H^k_2(\Omega)$ where $\alpha$ is a smooth representative of $c$, we choose 
$\bar\alpha\in C^\infty(\Lambda^kT^*M)$ a smooth extension of $\alpha$, then $d\bar\alpha$ is a closed
smooth form with support in $K$ and if $\iota\,:\,\partial K\rightarrow K$ is the inclusion we have
$\iota^*(d\bar\alpha)=0$. Some standard verifications show that 
$$b(c)=[d\bar\alpha]\in H^{k+1}(K,\partial K)$$ is well defined, that is it doesn't depend of the choice of
$\alpha\in c$ nor on the smooth extension of $\alpha$.

The inclusion map $j_\Omega\,:\, \Omega\rightarrow M$ induced a linear map (the restriction map) 
$$[j_\Omega^*]\,:\, H^k_2(M)\rightarrow H^k_2(\Omega).$$
\begin{lem}\label{exactabs}
We always have 
$$\ker b=\ima [j_\Omega^*].$$
\end{lem}
{\it Proof of lemma \ref{exactabs}.} First by construction we have $b\circ[j_\Omega^*]=0$ hence we only need to prove
that $$\ker b\subset\ima [j_\Omega^*].$$

Let $c\in \ker b$ and let $\alpha\in Z^k_2(\Omega)$ a smooth representative of $c$, we know that 
$\alpha$ has a smooth extension $\bar\alpha$ such that
$d\bar\alpha$ is zero in $H^{k+1}(K,\partial K)$. That is to say there is a smooth $k-$ form
$\beta\in C^\infty(\Lambda^kT^* K)$ such that 
$$d\bar\alpha=d\beta \mbox{ on } K \mbox{ and } \iota^*\beta=0$$
We claim that the $L^2$ form $\tilde  \alpha$ defined by 
$$\tilde \alpha=\left\{\begin{array}{ll} 
\bar\alpha-\beta& \mbox{ on } K\\
\alpha& \mbox{ on } \Omega\\
\end{array} \right.$$
is weakly closed.
As a matter of fact, we note
$\vec\nu\,:\, \partial K\rightarrow TM$ the unit normal vector field pointing into
$\Omega$ and $\iota\,:\, \partial K\rightarrow M$ the inclusion, then
 let $\varphi\in C^\infty_0(\Lambda^{k+1 }T^*M)$ with the Green's formula, we obtain
\begin{equation*}
\begin{split}
\int_M(\tilde\alpha,d^*\varphi) &=\int_K(\bar\alpha-\beta,d^*\varphi)
+\int_\Omega (\alpha,d^*\varphi)\\
&=-\int_{\partial K}(\iota^*(\bar\alpha-\beta), \inte_{\vec\nu}\varphi) d \sigma
+\int_{\partial K} (\iota^*\bar\alpha, \inte_{\vec\nu}\varphi) d \sigma\\
&=0.\\
\end{split}
\end{equation*}

We clearly have $j_\Omega^*\tilde \alpha=\alpha $, hence $c=[\alpha]\in \ima
[j_\Omega^*].$\hfill$\square$

Now because $H^{k+1}(K,\partial K)$ has finite dimension, we know that 
$$\dim\ima[j_\Omega^*]<\infty\Leftrightarrow \dim H^k_2(\Omega)<\infty.$$
Hence we get the implication :
$$\dim H^k_2(M)<\infty\Rightarrow \dim H^k_2(\Omega)<\infty.$$

To prove the reverse implication, we consider the reduced $L^2$ cohomology of $\Omega$ relative to the
boundary $\partial\Omega=\partial K$. We introduce 
\begin{equation*}
\begin{split}
Z^k_2(\Omega,\partial \Omega)&=\big\{\alpha\in L^2(\Lambda^kT^* \Omega), \mbox{ such that }
\forall \varphi \in C^\infty_0(\Lambda^kT^*  \overline{\Omega}),  \langle \alpha,d^*\varphi\rangle=0\big\}\\
&=\Big(d^* C^\infty_0(\Lambda^kT^*  \overline{\Omega})\Big)^\perp.\\
\end{split}
\end{equation*}

We remark here that the elements of $C^\infty_0(\Lambda^kT^*  \overline{\Omega})$ have compact support in 
$\overline{\Omega}$ in particular their support can touch the boundary; 
in fact a smooth $L^2$ closed $k-$form
belongs to $Z^k_2(\Omega,\partial \Omega)$ if and only its pull-back by $\iota$ is zero.
 This is a consequence of the integration by part formula
 $$\int_\Omega(\alpha,d^*\varphi) d\vol_g=\int_\Omega(d\alpha,\varphi)
 d\vol_g+\int_{\partial\Omega}(\iota^*\alpha,\inte_{\vec\nu}\varphi) d\sigma$$
 where $d\sigma$ is the Riemannian volume on $\partial \Omega$ induced by the metric $g$ and
 $\vec\nu\,:\, \partial\Omega\rightarrow T\Omega$ is the unit inward normal vector field.
 
We certainly have $dC^\infty_0(\Lambda^kT^* \Omega)\subset Z^k_2(\Omega,\partial \Omega)$ and we define
$$H^k_2(\Omega,\partial \Omega)
=\frac{Z^k_2(\Omega,\partial \Omega)}{\overline{dC^\infty_0(\Lambda^kT^* \Omega)}}\ .$$
These relative (reduced) $L^2$ cohomology space can be defined for every Riemannian manifold with boundary.
In fact, these relative (reduced) $L^2$ cohomology spaces  also have an interpretation in terms of harmonics
forms :

\begin{equation}\label{hodgerel}H^k_2(\Omega,\partial \Omega)\simeq \cH^k_{rel}(\Omega)\end{equation} where 
$$\cH^k_{rel}(\Omega)=\Big\{ \alpha\in L^2(\Lambda^kT^* \Omega), d\alpha=d^*\alpha=0 \mbox{ and }
\iota^*\alpha=0 \Big\}.$$ 

There is a natural map : the extension by zero map :
$$e\,:\, H^k_2(\Omega,\partial \Omega)\rightarrow H^k_2(M)\, .$$

When $\alpha\in L^2(\Lambda^kT^* \Omega)$ we define $e(\alpha)$ to be $\alpha$ on $\Omega$ and zero on $K$,
$e\,:\, L^2(\Lambda^kT^* \Omega)\rightarrow L^2(\Lambda^kT^* M)$ is clearly a bounded map moreover
$$e\Big(dC^\infty_0(\Lambda^kT^* \Omega)\Big)\subset dC^\infty_0(\Lambda^kT^* M),$$ hence by continuity of $e$
:
$e\Big(\overline{dC^\infty_0(\Lambda^kT^* \Omega)}\Big)\subset \overline{dC^\infty_0(\Lambda^kT^* M)}$.
Moreover when $\alpha\in Z_2^k(\Omega,\partial \Omega)$ then $e(\alpha)\in Z_ 2^k(M)$ :
because if $\varphi\in C^\infty_0(\Lambda^{k+1}T^* M)$ then
$$\langle e(\alpha),d^*\alpha\rangle=\int_\Omega (\alpha,d^*(j_\Omega^*\varphi))d\vol_g$$
but $(j_\Omega^*\varphi)\in C^\infty_0(\Lambda^kT^*  \overline{\Omega})$ hence this integral is zero by definition 
of $Z_2^k(\Omega,\partial \Omega)$.

Let $j_K\,:\, K\rightarrow M$ the inclusion map, it induces as before a linear map
$$[j_K^*]\,:\, H^k_2(M)\rightarrow H^k_2(K) \simeq H^k(K).$$
We always have $[j_K^*]\circ e=0$ hence 
$$\ima e\subset \ker [j_K^*]\ .$$ In fact as before 
\begin{lem}\label{exactrel}We have the equality
$$\ima e= \ker [j_K^*].$$\end{lem}
{\it Proof of the lemma \ref{exactrel}.}
If $c\in \ker [j_K^*]$ and $\alpha\in L^2(\Lambda^kT^* M)$ is a smooth representative of $c$ (for instance $\alpha$ is $L^2$ and
harmonic). By definition we know that there is $\beta\in C^\infty(\Lambda^{k-1}T^*K)$ such that
$$j_K^*\alpha=d\beta$$
Now consider $\bar\beta$ any smooth extension of $\beta$ with compact support.
We clearly have $$[\alpha-d\bar\beta]=[\alpha]=c\mbox{ in } H^k_2(M).$$
Moreover by construction 
$$\alpha-d\bar\beta=e(j_\Omega^*(\alpha-d\bar\beta))$$
If we verify that $j_\Omega^*(\alpha-d\bar\beta)\in Z^k_2(\Omega,\partial\Omega)$ we have finish the proof of the
equality the lemma \ref{exactrel}. In fact, this verification is straightforward.
Let $\varphi\in C^\infty_0(\Lambda^kT^*  \overline{\Omega})$ and consider $\bar\varphi$ any smooth extension of
$\varphi$ :
$$\int_\Omega(j_\Omega^*(\alpha-d\bar\beta),d^*\varphi)d\vol_g
=\int_M(\alpha-d\bar\beta,d^*\bar\varphi)d\vol_g=0\ .$$
\hfill$\square$

Again since $H^k(K)$ have finite dimension, the kernel of $[j_K^*]$ have finite codimension in $H^k_2(M)$
and we have obtain the implication 
$$\dim H^k(\Omega,\partial \Omega)<\infty\Rightarrow \dim H^k_2(M)<\infty.$$

In order to conclude, we use the Hodge star operator; because our manifold is oriented, the Hodge star
operator is an isometry which exchanges $k$ forms and $(n-k)$ forms :
$$\star\,:\, \Lambda^kT^*_xM\rightarrow \Lambda^{n-k}T^*_xM.$$
This operator satisfies the following properties
 $$\star\circ\star=\pm \Id$$
 where the sign depends on the degree.
 Moreover we have :
 $d^*=\pm \star d\star$. 
 Hence the Hodge star operator maps the space $\cH^k(M)$ to $\cH^{n-k}(M)$. Moreover
 it is also a good exercise to check that in our setting :
 $$\iota^*(\star\alpha)=\pm \inte_{\vec\nu}\alpha.$$
 Hence the Hodge star operator maps the space $\cH_{abs}^k(\Omega)$ to $\cH_{rel}^{n-k}(\Omega)$.
 
But our last result says that 
 $$\dim \cH^k_{rel}(\Omega)<\infty\Rightarrow \dim \cH^k_2(M)<\infty.$$
 Hence using the Hodge star operator we get
 $$\dim \cH^{n-k}_{abs}(\Omega)<\infty\Rightarrow \dim \cH^{n-k}_2(M)<\infty.$$
 That is $\dim H_2^{n-k}(\Omega)<\infty\Rightarrow \dim H^{n-k}_2(M)<\infty.$
  It is now clear that we have prove Lott's results.\hfill$\square$

\begin{rems}
\begin{enumeroman}
\item The reader can verify that using forms with coefficients in the orientation bundle, we can removed the
orientability condition.
\item These two properties are in fact a heritage of the following two exact sequences for the de Rham cohomology:
$$...\rightarrow H^k(K,\partial K)\rightarrow H^k(M)\rightarrow H^k(\Omega)\rightarrow H^{k+1}(K,\partial
K)\rightarrow...$$
$$...\rightarrow H^k(\Omega,\partial \Omega)\rightarrow H^k(M)\rightarrow H^k(K)\rightarrow 
H^{k+1}(\Omega,\partial \Omega)\rightarrow...$$
\end{enumeroman}
\end{rems}
\subsection{Some bibliographical hints}

We recommend the reading of the classical book of G. de Rham \cite{dR} or
W. Hodge \cite{Hodge}. For a modern treatment of the de Rham's theorem a very good reference is the book
 of Bott-Tu \cite{BT} . The Hodge theorem for compact manifold
with boundary have been proved by P.E. Conner \cite{conn}. Other proofs of the 
Hodge-de Rham theorem can be found in other classical book (for instance
in the book of Griffith-Harris (chapter 0 section 6 in \cite{GH}) or in the book
of  M. Taylor (chap. 5 in \cite{T}). We should also mentioned a sheaf theoretical proof of the 
Hodge-de Rham theorem by N. Telemann \cite{Te} . About the general feature on non compact
manifold, you can read the paper of J. Dodziuk \cite{Dod},  the now classical paper of J. Cheeger 
\cite{Che}, the first section of the paper by J. Lott \cite{lottL2} or look at
\cite{carrroma} and also read the beautiful paper of M. Anderson \cite{An}. The paper of J. Br\"uning and M.
Lesch deals with an abstract approach about the
identification between the space of $L^2$ harmonic form and $L^2$ cohomology \cite{BL}.
For a first approach on the $L^2$ cohomology of symmetric space,
 the paper of S. Zucker \cite{Z1} is very nice, the
first parts of the survey of W. Casselman is also instructive \cite{cass}.
\section{Harmonics $L^2$ $1-$ forms}

In this second lecture, we survey some results between the space of 
harmonics $L^2$ $1-$ forms, cohomology with compact support and the geometry of ends.

\subsection{Ends}
\subsubsection{Definitions}
When $U\subset M$ is a subset of an open manifold $M$ we say that $U$ is {\it bounded} if $\overline{U}$ is a
compact subset of $M$, when $U$ is not bounded we say that $U$ is {\it unbounded}. When $g$ is complete
Riemannian metric on $M$, then $U\subset M$ is bounded if and only if there is some $R>0$ and $o\in M$ such
that $U\subset B(o,R)$.

Let $M$ be a smooth manifold, we say that $M$ has only one end if for any compact subset $K\subset M$,
$M\setminus K$ has only one unbounded connected component. 

For instance, when $n\ge 2$, the Euclidean space $\R^n$ has only one end.

More generaly, we say that $M$ has  $k$ ends (where $k\in \N$) if there is a compact set $K_0\subset M$ such
that for every compact set $K\subset M$ containing $K_0$, $M\setminus K$ has exactly 
$k$ unbounded connected components. 

 A compact manifold is a manifold with zero end.
 For instance, $\R$ or $\R\times\bS^1$ have two ends. The topology at infinity of manifold with 
  only one end can be very complicated. 
  
  \subsubsection{Number of ends and cohomology}
   Here $M$ is a smooth non compact connected manifold.
   \begin{lem} If $M$ has at least two ends, then $$H^1_0(M)\not=\{0\}.$$
   In fact, if $M$ has at least $k$ ends then
   $$\dim H^1_0(M)\ge k-1.$$
   \end{lem}
\proof If $M$ has more than two ends, then we can find a compact set $K\subset M$ such that
$$M\setminus K=U_-\cup U_+$$ where
$U_-,U_+$ are unbounded and $U_-\cap U_+=\emptyset.$
\begin{center}\fbox{\includegraphics[height=6cm]{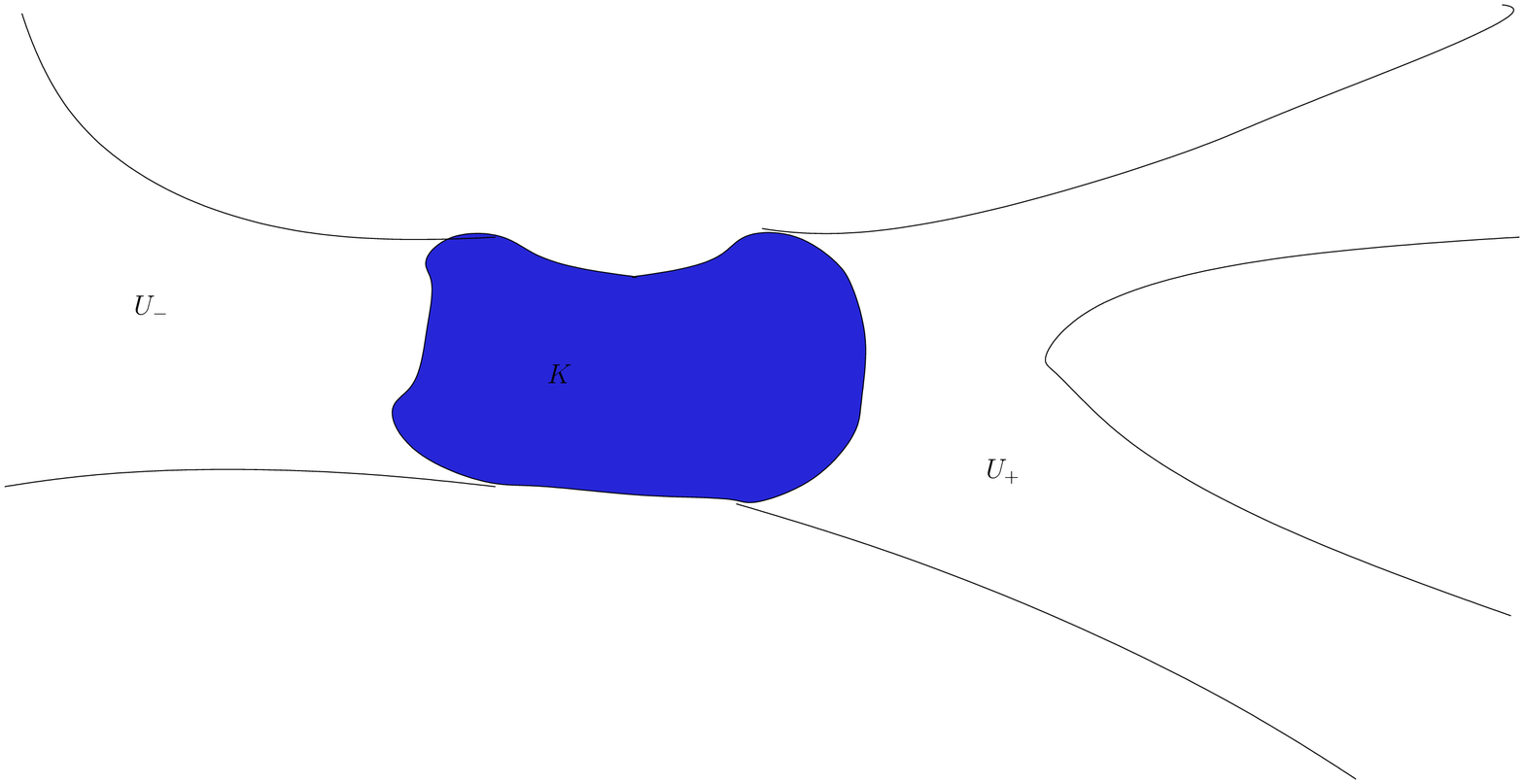}  }
\end{center} Let $u\in C^\infty(M)$ such that 
$u=\pm 1$ on $U_\pm$, then clearly $\alpha=du$ is a closed $1-$form with compact support.
If the cohomology class of $\alpha$ in $H^1_0(M)$ is trivial, then we
 find some $f\in C^\infty_0(M)$ such that $du=\alpha=df$. But $M$ is connected hence
we have a constant $c$ such that $u=f+c$. We look at this equation outside the support of $f$ and
on $U_\pm$ we find that $c=\pm 1.$
\endproof

In fact there is a weak reciprocal to this result
(see proposition 5.2 in \cite{carrped}) :
\begin{prop} \label{cohomZ}
If $M^n$ is an open manifold having one end, and if every twofold 
normal covering of
$M$ has also one end, then $$H^1_0(M,\Z)=\{0\}.$$
In particular, $H^1_0(M)=\{0\}$ and if furthermore
$M$ is orientable, then
$$H_{n-1}(M,\Z)=\{0\}.$$
\end{prop}

\subsection{$H^1_0(M)$ versus $H^1_2(M)$}
From now we assume that $(M,g)$ is a complete Riemannian manifold.
\subsubsection{An easy case}
\begin{lem}Assume that $M$ has only one end, then 
$$\{0\} \rightarrow H^1_0(M)\rightarrow H^1_2(M).$$
\end{lem}
\proof Let $\alpha\in C^\infty_0(T^*M)$ be a closed $1-$form which is zero in
$H^1_2(M).$ By the result (\ref{dRham}), we know that there is a smooth function $f\in C^\infty(M)$ such that
$$\alpha=df.$$
But $M\setminus \supp\alpha $ has only one unbounded connected component $U$, hence
on $U$ $df=0$ hence there is a constant $c$ such that $f=c$ on $U$. Now by construction the function
$f-c$ has compact support and $\alpha=d(f-c).$ Hence $\alpha$ is zero in $H^1_0(M)$ \endproof
In fact the main purpose of this lecture is to go from geometry to topology : we want to find
geometrical conditions insuring that this map is injective. 

\subsubsection{Condition involving the spectrum of the Laplacian}
\begin{prop}\label{spectr} Assume that all ends of $M$ have infinite volume 
\footnote{that is ouside every compact subset of $M$, all unbounded connected
components have infinite volume.} and assume that there is a $\lambda>0$ such that 
\begin{equation}
\label{gap0}\forall f\in C^\infty_0(M),\ \lambda\int_M f^2 d\vol_g\le \int_M |df|^2 d\vol_g.\end{equation}
Then $$\{0\} \rightarrow H^1_0(M)\rightarrow H^1_2(M).$$
\end{prop}
\proof Let $[\alpha]\in H^1_0(M)$ is mapped to 
zero in $H_ 2^1(M)$. Hence
there is a sequence $(f_k)$ of smooth functions with compact support on $M$
such that $\alpha=L^2-\lim df_k$. Since we have the inequality
$$\|df_k-df_l\|_{L^2}^2\ge \lambda \|f_k-f_l\|_{L^2}^2,$$
and since $\lambda>0$, we conclude that this sequence $(f_k)$
converges to some $f\in L^2$, so that
$\alpha=df$.
But $\alpha$ has compact support, hence $f$ is locally constant 
outside the compact
set $\supp(\alpha)$. Since all unbounded connected components of $M\priv
\supp(\alpha)$ have infinite volume and since $f\in L^2$, we see that $f$ has
compact support, hence $[\alpha]=[df]=0$ in $H^1_0(M)$.\endproof

Before giving some comments on the hypothesis on this proposition, let's give a consequence of this
injectivity :
\begin{prop} Assume that $$\{0\} \rightarrow H^1_0(M)\rightarrow H^1_2(M),$$
and that $M$ has at least $k$ ends then
\begin{equation}\label{dimHD}\dim\big\{h\in C^\infty(M), \Delta h=0\mbox{ and }dh\in L^2\big\}\ge k.\end{equation}
\end{prop}
A general formula for the dimension of the space of bounded harmonic function with $L^2$ gradient
has been obtained by A. Grigor'yan (\cite{Gri0}).
\proof Assume that $k=2$ (the other cases are similar), there is a compact set $K\subset M$ such that
$$M\setminus K=U_-\cup U_+$$ with 
$U_-,U_+$ unbounded and $U_-\cap U_+=\emptyset.$ Let $u\in C^\infty(M)$ such that 
$u=\pm 1$ on $U_\pm$ we know that $\alpha=du$ is not zero in cohomology with compact support , hence we know
that there is a non zero harmonic  $L^2$ $1-$form $\eta$  such that
$[\alpha]=[\eta]$ in $H^1_2(M)$; in particular with (\ref{dRham}), we can find $v\in C^\infty(M)$
such that 
$$du=\alpha=\eta+dv$$
that is if $h=u-v$ then 
$dh=\eta\in L^2$ and 
$\Delta h=d^*d(u-v)=d^*\eta=0$. And because $\eta\not=0$ we know that $h$ is not the constant function. The
linear span of $h$ and of the constant function $1$ is
 of dimension $2$ and is include in 
 $$ \dim\big\{h\in C^\infty(M), \Delta h=0\mbox{ and }dh\in L^2\big\}.$$\endproof
 
 \begin{rem}\label{alter}
 In the setting of the proposition (\ref{spectr}), we can give a direct and more classical
 proof of this inequality : we assume that $(M,g)$ satisfies that the assumption made in proposition
 (\ref{spectr}), we are going to show the inequality (\ref{dimHD}).
Again we assume that $k=2$ : there is a compact set $K\subset M$ such that
$$M\setminus K=U_-\cup U_+$$ with 
$U_-,U_+$ unbounded and $U_-\cap U_+=\emptyset.$
Let $o\in M$ be a fixed point and for large $k$ we have
$K\subset B(o,k)$. We consider the solution of the Dirichlet problem :
$$\left\{\begin{array}{ll}
\Delta u_k=0 &\mbox{ on } B(o,k)\\
u_k=\pm 1 &\mbox{ on } \partial B(o,k)\cap U_\pm\\
\end{array}\right.$$
Then $u_k$ is the minimizer of the functional
$$u \mapsto \int_{B(o,k)} |du|^2$$
amongst all functions in $W^{1,2}(B(o,k))$ (functions in $L^2(B(o,k))$ whose derivatives are also in $L^2$)
such that $u=\pm 1 $ on $ U_\pm\cap\partial B(o,k)$\,). We extend $u_k$ to all
$M$ by setting $u_k=\pm 1$ on $ U_\pm\setminus B(o,k)$.
Then this extension (also denote by $u_k$) is the minimizer of the functional
$$u \mapsto \int_{M} |du|^2$$ amongst all functions in $W^{1,2}_{loc}$ 
(functions in $L^2_{loc}$ whose derivatives are also in $L^2_{loc}$)
such that $u=\pm 1 $ on $U_\pm\setminus B(o,k).$
Hence we always have when $k<l$ :
$$\int_{M} |du_l|^2\le\int_{M} |du_k|^2.$$
Moreover by the maximum principle, we always have
$$-1\le u_k\le 1.$$
After extraction of a subsequence we can assume that
uniformly on compact set 
$$\lim_{k \to\infty} u_k=u.$$
The function $u$ is then a harmonic function whose value are in $[-1,1]$ and it satisfies
$$\int_{M} |du|^2\le\int_{M} |du_k|^2<\infty.$$
We must show that this $u$ is not a constant function. We apply our estimate to the function $u_k-u_l$ where
$l\ge k$, then we get 
$$\lambda\int_M |u_k-u_l|^2\le \int_M |du_k-du_l|^2\le 4\int_M |du_k|^2$$
In particular if we let $l\to\infty$, we find that
$$u_k-u\in L^2$$ that is
$$\int_{U_\pm} |u-(\pm1)|^2 <\infty.$$
But by hypothesis the volume of $U_\pm$ are infinite hence $u$ cannot be the constant function.
\end{rem}
We know make some comments on the assumption of the proposition (\ref{spectr}) :

\noindent{\it On the first condition.} We have the following useful criterion :
\begin{lem}\label{volumeend} Assume that there is $v>0$ and $\varepsilon>0$ such that for all 
$x\in M$ 
$$\vol B(x,\varepsilon)\ge v$$ then all ends of $M$ have infinite volume.
\end{lem}
\proof We fix a based point $o\in M$ and let $K\subset M$ be a compact set
 and $U$ an unbounded connected component of $M\setminus K$. Let
$R>0$ large enough so that $K\subset B(o,R)$ (recall that we have assumed here that the Riemannian manifold
$(M,g)$ is complete). For $k\in \N$ we choose $x\in \partial B(o,R+(2k+1)\varepsilon)\cap U$ and we consider
$\gamma\,:\,[0,R+(2k+1)\varepsilon]\rightarrow M$ a minimizing geodesic from $o$ to $x$ (parametrize by arc-
length). Necessary for all $t\in ]R,R+(2k+1)\varepsilon]$ we have $\gamma(t)\in U$ ; moreover for 
$l=0,1,..., k$,
the open geodesic ball $B(\gamma(R+(2l+1)\varepsilon), \varepsilon)$ are in $U$ and disjoint
hence
$$\vol U\ge \sum_{l=0}^k \vol B\big(\gamma(R+(2l+1)\varepsilon),\varepsilon \big)\ge (k+1)v. $$
\begin{center}\fbox{\includegraphics[height=6cm]{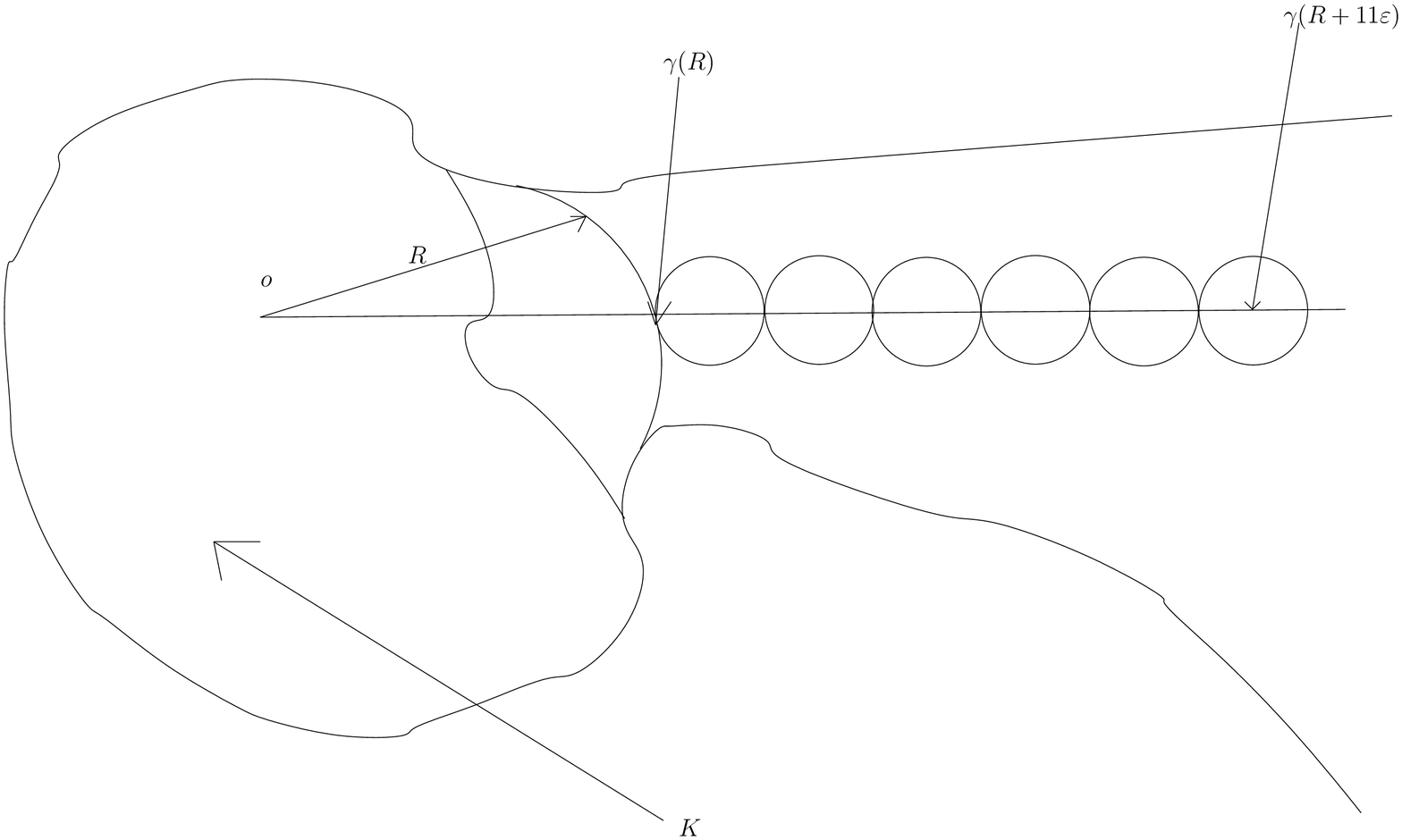}  }
\end{center}
\endproof
For instance, when the injectivity radius of $(M,g)$ is positive, C. Croke has shown in (\cite{croke}) that for all
$x\in M$ and all $r\le \inj (M)$ then
$$\vol B(x,r)\ge C_n r^n.$$ Hence  a complete Riemannian manifold with positive injectivity radius has all
its ends with infinite volume.

\noindent{\it About the second condition} The condition (\ref{gap0}) is clearly equivalent to
$$\lambda_0(M,g)
:=\inf_{f\in C^\infty_0(M)}\left\{\frac{\int_M |df|^2 d\vol_g}{\int_M |f|^2 d\vol_g}\right\}>0.$$
This condition is linked with the spectrum of the Laplacian on functions. In fact $\lambda_0(M,g)$ is the
bottom of the spectrum of the operator $d^*d=\Delta$ acting on functions. 
Because the Riemannian manifold
$(M,g)$ is complete, the operator 
$$d^*d=\Delta\,:\, C^\infty_0(M)\rightarrow L^2(M,d\vol_g)$$ has a unique self adjoint extension with domain :
$$\cD(\Delta)=\{f\in L^2(M), \Delta f\in L^2(M)\}=\{ f\in L^2(M), df\in L^2(M)\mbox{ and } d^*df\in L^2.\}$$
Hence $f\in \cD(\Delta)$ if and only if $f\in \cD^0(d)$ and there is $C\in \R$ such that 
$$ \forall\varphi\in C^\infty_0(M)\ ,\ |\langle df,d\varphi\rangle|\le C \|\varphi\|_{L^2}.$$

The spectrum of $\Delta_g$ is a closed subspace of $[0,\infty[$ and the condition $\lambda_0(M,g)>0$ is
equivalent to zero not being in the spectrum of $\Delta$. 
Also this condition depends only of the geometry at infinity, in fact we have that $\lambda_0(M,g)>0$
if and only if there is
 $K\subset M$ a compact
subset such that $\lambda_0(M\setminus K,g)>0$ and $\vol(M)=\infty$.. 
From the proof of the proposition (\ref{spectr}) it is clear that we have the following result :
\begin{prop}
Assume that $(M,g)$ is a complete Riemannian manifold such that 
$\lambda_0(M,g)>0$, then if we introduce 
$$H^1_v(M)=\frac{\{\alpha\in C^\infty_0(T^*M), \mbox{such that }  d\alpha=0\}}{\{df, f\in C^\infty(M),
\mbox{such that }df\in C^\infty_0(T^*M) \mbox{ and }\vol(\supp f)<\infty\}}$$ then
$$\{0\} \rightarrow H^1_v(M)\rightarrow H^1_2(M).$$
\end{prop}

{\it Optimality of the results :} 
The real line $(\R, (dt)^2)$ has clearly all his ends with infinite volume but 
$\lambda_0(\R, dt^2)=0$ ; as a matter of fact  if $u\in C^\infty_0(\R)$ is not the zero function, then
for $u_n(t)=u(t/n)$ we have
$$\int_\R |u_n'(t)|^2 dt =\frac{1}{n}\int_\R |u'(t)|^2 dt$$
and 
$$\int_\R |u_n(t)|^2 dt =n \int_\R |u'(t)|^2 dt.$$

But the first group of cohomology with compact support of $ \R$ has dimension $1$ where as
$$\cH^1(\R,dt^2)=\{fdt, f\in L^2\mbox{ and } f'=0\}=\{0\}.$$

We consider now the manifold $\Sigma=\R\times\bS^1$ endowed with the warped product metric :
$$g=(dt)^2+e^{2t}(d\theta)^2.$$
$M$ have two ends one with finite volume and the other of infinite volume. Again the 
first group of cohomology with compact support of $\Sigma$ has dimension $1$. 
Moreover $$\lambda_0(\Sigma)\ge \frac14\ ;$$
as a matter of fact let $f\in C^\infty_0( \Sigma)$
 $$\int_\Sigma |df|^2 d\vol_g
 \ge \int_\Sigma \left|\frac{\partial f}{\partial t}(t,\theta)\right|^2 e^{t}dtd\theta$$
 Let $f=e^{-t/2} v$ then 
\begin{equation*}
\begin{split}
\int_\Sigma \left|\frac{\partial f}{\partial t}(t,\theta)\right|^2 e^{t}dtd\theta&=
\int_\Sigma \left|\frac{\partial v}{\partial t}-\frac12 v\right|^2 dtd\theta \\
&=\int_{\R\times\bS^1} \Big[
\left|\frac{\partial v}{\partial t}\right|^2+\frac14 \left|v\right|^2+v\frac{\partial v}{\partial t}
\Big] dtd\theta\\
&=\int_{\R\times\bS^1} \Big[
\left|\frac{\partial v}{\partial t}\right|^2+\frac14 \left|v\right|^2+\frac12\frac{\partial v^2}{\partial t}
\Big] dtd\theta\\
&=\int_{\R\times\bS^1} \Big[
\left|\frac{\partial v}{\partial t}\right|^2+\frac14 \left|v\right|^2
\Big] dtd\theta\\
&\ge \frac14\int_{\R\times\bS^1} \left|v\right|^2
 dtd\theta\\
 &=\frac14\int_\Sigma |f|^2 d\vol_g.\\
\end{split}
\end{equation*}
In fact we can show (using an appropriate test function
$u_n(t,\theta)=\chi_n(t) e^{-t/2}$) that $\lambda_0(\Sigma)=\frac14\ $.

In fact we can show that $\cH^1(\Sigma)$ has infinite dimension (see the first part of theorem
\ref{surface}), however if 
the conclusion of the proposition \ref{spectr}) were true for $(\Sigma,g)$ 
then by the alternative proof of (\ref{dimHD}) in (\ref{alter}), we would find
 a non constant harmonic function $h$ with $L^2$ gradient.
If we look at the construction of this function, it is not hard to check that
$h$ will only depends on the first variable that is
$$h(t,\theta)= f(t)$$ where $f$ solves the O.D.E. 
$$(e^tf')'=0 \,;$$ hence there is constant $A$ and $B$ such that
$$h(t,\theta)=Ae^{-t} +B.$$
The fact that $dh\in L^2$ implies that 
$$\int_\R (f')^2(t) e^t dt<\infty$$ that is
$A=0$ and $h$ is the constant function.

This simple proposition \ref{spectr} or variant of it has been used frequently, I will
mention only two results, the first one is the following very beautifully result of P. Li and J. Wang
(\cite{LW1}):
(see also the other articles of P. Li and J. Wang for other related results \cite{LW2,LW4})
\begin{thm} Assume that $(M,g)$ is a complete Riemannian manifold of dimension $n>2$ with 
$$\lambda_0(M)\ge (n-2)$$ assume moreover that its Ricci curvature satisfies the lower bound
$$\ricci_g\ge -(n-1) g$$
then either $(M,g)$ has only one end of infinite volume either
$(M,g)$ has two ends with infinite volume and  is isometric to the warped product 
$\R\times N$ endowed with the metric $$(dt)^2+\cosh^2(t) h$$ where $(N,h)$ is a compact Riemannian manifold
with $\ricci_g\ge -(n-2) g$
\end{thm}

The second concerns the locally symmetric space (\cite{carrped}) :
\begin{thm} \label{endT}
Let $G/K$ be a symmetric space without any compact factor and without 
any factor
isometric to a real or complex hyperbolic space. Assume that 
$\Gamma\subset G$ is a torsion-free,
discrete subgroup of $G$ such that $\Gamma\backslash G/K$ is non 
compact and that all ends  of
$\Gamma\backslash G/K$ have infinite volume. Then
$\Gamma\backslash G/K$ has only one end, and 
$$H_{m-1}(\Gamma\backslash G/K,\Z)=\{0\},$$
where $m=\dim(G/K)$.
\end{thm}

\subsubsection{Condition involving a Sobolev inequality }
\begin{prop}\label{sobolev} Assume that $(M,g)$ is a complete manifold that satisfies for a $\nu >2$ and $\mu>0$ the Sobolev
inequality :
\begin{equation}
\label{sob}\forall f\in C^\infty_0(M),\ \ \mu\left(\int_M f^{\frac{2\nu}{\nu-2}}
d\vol_g\right)^{1-\frac{2}{\nu}}\le \int_M |df|^2 d\vol_g.\end{equation}
Then $$\{0\} \rightarrow H^1_0(M)\rightarrow H^1_2(M).$$
\end{prop}
The proof follows essentially the same path : if $[\alpha]\in H^1_0(M)$ is map to zero in $H^1_2(M)$,
then we find $f\in L^{\frac{2\nu}{\nu-2}}$ such that $$\alpha=df.$$
Hence $f$ is locally constant outside the support of $\alpha$. 
But according to (proposition 2.4 in \cite{car-smf}), we know that the Sobolev inequality (\ref{sob}) implies a uniform lower
bound on the volume of geodesic balls : 
$$\forall x\in M, \forall r\ge 0\ :\ \vol B(x,r)\ge C(\nu) \left(\mu r^{2}\right)^{\nu/2}$$
for some explicit constant $C(\nu)>0$ depending only on $\nu$.
Hence by (\ref{volumeend}), we know that such an estimate implies that all the unbounded connected
components of $M\setminus\supp\alpha$ have infinite volume. Hence $f$ has compact support and 
$[\alpha]=0 $ in $H^1_0(M)$.

Examples of manifolds satisfying the Sobolev inequality :
\begin{enumeroman}
\item The Euclidean space $\R^n$ $(n\ge 3)$ satisfies the Sobolev inequality for $\nu=n$.
\item If $M^n\subset \R^N$ is a minimal submanifold then $M^n$ with the induced metric satisfies
the Sobolev inequality for $\nu=n$ \cite{MS,HS}.
\item A Cartan-Hadamard manifold (a simply connected non positively curved complete Riemannian manifold) of
dimension $n\ge 3$ satisfies
the Sobolev inequality for $\nu=n$.
\item The hyperbolic space of dimension $n$, satisfies the Sobolev inequality for any
 $\nu\in[n,\infty [\cap
]2,\infty[.$
\end{enumeroman}

In fact the validity of the Sobolev inequality depends only of the geometry at infinity : 
according to \cite{carrduke}, if $K\subset M$ is a compact subset of $M$ then the Sobolev inequality 
(\ref{sob}) holds for $M$ if and only if it holds for $M\setminus K$:
$$\forall f\in C^\infty_0(M\setminus K),\ \ \mu\left(\int_{M\setminus K} f^{\frac{2\nu}{\nu-2}}
d\vol_g\right)^{1-\frac{2}{\nu}}\le \int_{M\setminus K} |df|^2 d\vol_g.$$

This proposition (\ref{sobolev}) has the following beautiful application due to Cao-Shen-Zhu \cite{CSZ} (see also \cite{LW3})
however these authors haven't prove the vanishing of the first group of cohomology with compact support :
\begin{cor} If $M^n\subset \R^{n+1}$ is a stable complete minimal hypersurface then $M$ has only one end,
moreover $H^1_0(M)=\{0\}.$
\end{cor}
\proof We only need to show that $\cH^1(M)=\{0\}$. Let $\alpha\in\cH^1(M)$ then it satisfies the Bochner
identity 
$$\int_M \big[\left|\nabla\alpha\right|^2+\ricci (\alpha,\alpha)\big] d\vol_g=0.$$
However the Gauss equations imply that 
$$\left|\ricci (\alpha,\alpha)\right|\le |A|^2 \left|\alpha\right|^2,$$
where $A$ is the second fundamental form of the hypersurface $M\subset \R^{n+1}$. 
The stability condition says that the second variation  of the area is non negative that is
$$\forall f\in C^\infty_0(M),\ J(f):=\int_M \big[|df|^2-|A|^2\,|f|^2\big] d\vol_g\ge 0.$$
But the refined Kato inequality  (due in this case to S.T. Yau) shows that
$$\left|\nabla\alpha\right|^2\le \frac{n}{n-1} \big|d|\alpha|\big|^2.$$

Hence we get 
\begin{equation*}
\begin{split}
0&=\int_M \big[\left|\nabla\alpha\right|^2+\ricci (\alpha,\alpha)\big] d\vol_g\\
&\ge \int_M \left[\frac{n}{n-1} \left|d|\alpha|\right|^2+\ricci (\alpha,\alpha)\right] d\vol_g\\
&\ge \int_M \left[ \frac{1}{n-1} \left| d|\alpha| \right|^2 
+\left|d|\alpha|\right|^2-|A|^2 \left|\alpha\right|^2 \right] d\vol_g\\
&\ge J(|\alpha|)+\frac{1}{n-1} \int_M\left|d|\alpha|\right|^2  d\vol_g.\\
&\ge \frac{1}{n-1} \int_M\left|d|\alpha|\right|^2  d\vol_g.\\
\end{split}
\end{equation*}
Hence $\alpha$ has constant length, but because of the Sobolev inequality, 
the volume of $(M,g)$ is infinite and
$\alpha=0$.
\endproof

\subsubsection{What is behind the injectivity of the map $H^1_0(M)\rightarrow H^1_2(M)$ }
In fact there is a general notion from potential theory which is related to the injectivity of this map
(see the survey of A. Ancona \cite{Ancona} or  \cite{LT}).
\begin{defi}\label{para} Let $E \subset M$ be an open connected set with smooth compact boundary, the following properties are
equivalent:
\begin{enumeroman}
\item there is a positive super-harmonic function : $s\,:\,\overline{E}\rightarrow \R_+^*$ with 
$$\liminf_{x\to\infty} s(x)=\inf_{x\in E} s(x)=0$$ and
$$\inf_{x\in\partial E} s(x)\ge 1.$$
\item There is a positive harmonic function : $h\,:\,\overline{E}\rightarrow \R_+^*$ with 
$$\liminf_{x\to\infty} h(x)=\inf_{x\in E} h(x)=0$$ and $h=1$ on $\partial E$.
Moreover $dh\in L^2$.
\item The capacity of $E$ is positive :$$\capa (E)=
\inf\left\{\int_E |dv|^2,\ v\in C^\infty_0(\overline{E})\ \mbox{and } v\ge 1 \mbox{ on }\partial
E\right\}>0.$$
\item For any $U\subset E$ bounded open subset of $E$ there is a constant $C_U>0$ such that 
$$\forall f\in  C^\infty_0(\overline{E}), C_U\int_U f^2\le\int_E |df|^2.$$
\item For some $U\subset E$ bounded open subset of $E$ there is a constant $C_U>0$ such that 
$$\forall f\in  C^\infty_0(\overline{E}), C_U\int_U f^2\le\int_E |df|^2.$$
\end{enumeroman}
When one of these properties holds we say that $E$ is non parabolic and if one of these properties fails we say
that $E$ is parabolic.
\end{defi}
\proof We clearly have $ii)\Rightarrow i)$ and $iv)\Rightarrow v)$.
We first prove that $i)\Rightarrow iv)$.
Let $\vec \nu\,:\, \partial E\rightarrow TE$ be the unit inward normal vector field along $\partial E$. 
When $v \in C^\infty_0(\overline{E})$ we set $\varphi=v/\sqrt{s}$ so that 
$$\int_E |dv|^2=\int_E s|d \varphi|^2+\int_E \varphi \langle ds, d\varphi\rangle+\int_E \frac{|ds|^2}{4s}
\varphi^2.$$
But 
\begin{equation*}
\begin{split}
2\int_E \varphi \langle ds, d\varphi\rangle&=\int_E \langle ds, d\varphi^2\rangle\\
&=\int_E(\Delta s) \varphi^2-\int_{\partial E} \varphi^2ds(\vec\nu) .\\
\end{split}
\end{equation*}
But $s$ is assumed to be super-harmonic hence $\Delta s\ge 0$ and $ds(\vec\nu)\le 0$ along $\partial E$, hence
we finally obtain the lower bound 
$$\int_E |dv|^2\ge \int_E \frac{|ds|^2}{4s}\varphi^2=\frac14 \int_E|d\log s|^2 v^2.$$
By assumption, $s$ is not constant hence we can find an open bounded set $U\subset E$ and $\varepsilon>0$
such that 
$|d\log s|>\varepsilon$ on $U$ and we obtain that for all $v \in C^\infty_0(\overline{E})$
$$\int_E |dv|^2\ge \frac{\varepsilon^2}{4}\int_U v^2.$$

We prove now that $v)\Rightarrow iv).$
Let $U\subset E$ such as in iv). For $o\in M$ a fixed point and
$R$ such that $U\cup \partial E\subset B(o,R)$, we will prove that there is a constant $C_R>0$
such that 
$$\forall f\in  C^\infty_0(\overline{E}),\ C_R\int_{ B(o,R)\cap E} f^2\le\int_E |df|^2.$$
Let $V\subset U$ be an non empty open set with $\bar V\subset U$ and 
let $\rho\in C^\infty(\overline{E})$ such that $\supp \rho\subset U$ and $0\le \rho\le 1$ and $\rho=1$ in $V$. Then
for $f\in  C^\infty_0(\overline{E})$ we have :
\begin{equation*}
\begin{split}
\int_{B(o,R)\cap E}f^2&\le  2\int_{B(o,R)\cap E}(\rho f) ^2+2\int_{B(o,R)\cap E}((1-\rho) f) ^2\\
&\le  2\int_{U} f^2+2\int_{B(o,R)\cap E}((1-\rho) f) ^2\\
&\le \frac{2}{C_U}\int_{ E}\left|df\right| ^2+
\frac{2}{\lambda}\int_{B(o,R)\cap E}\left|d((1-\rho)f)\right|^2\\
\end{split}
\end{equation*}
Where $\lambda>0$ is the first eigenvalue of the Laplacian on functions on $(B(o,R)\cap E)\setminus V$ for the Dirichlet
boundary condition on $\partial V$.
But 
\begin{equation*}
\begin{split}\int_{B(o,R)\cap E}\left|d((1-\rho)f)\right|^2&\le
2\int_{B(o,R)\cap E}\left|df\right|^2+2\|d\rho\|^2_{L^\infty}\int_{U}\left|f\right|^2 \\
&\le \left(2+\frac{2\|d\rho\|^2_{L^\infty}}{C_U}\right)\int_{ E}\left|df\right|^2.\\
\end{split}
\end{equation*}
Hence the result for $$C_R=\left(\frac{2}{C_U}+\frac{4}{\lambda}\left(1+\frac{\|d\rho\|_{L^\infty}^2}{C_U} \right)\right)^{-1}.$$

Now we consider the implication $iii)\Rightarrow ii)$.
We introduce 
$$C(R)=\inf\int_{E\cap B(o,R)} |dv|^2$$
where the infimum runs over all functions $v\in C^\infty(\overline{E}\cap B(o,R) )$
such that $ v\ge 1$ on $\partial E$ and $v=0$ on $\partial B(o,R)\cap E$ ;
where $o\in M$ is a fixed point and where $R>0$ is chosen large enough so that $\partial E\subset B(o,R).$
We have $C(R)>0$ and we have assumed that $C(\infty)=\inf_R C(R)>0$. Each $C(R)$ is realized by the
harmonic function $h_R$ such that 
$h_R=1$ on $\partial E$ and $h_R=0$ on $\partial B(o,R)\cap E$.
We extend $h_R$ by zero on $E\setminus B(o,R)$, an application of the maximum principle implies that when 
$R\ge R'$ then $h_{R'}\le h_R$. We let
$h(x)=\sup_{R} h_R(x)$. On compact subset of $\overline{E}$, $h_R$ converge to $h$ in the smooth topology
moreover $h$ is a harmonic function with $0\le h\le 1$. The Green formula shows that
$$C(R)=-\int_{\partial E} dh_R(\vec\nu),$$
hence when $R\to\infty$ we obtain
$$C(\infty)=-\int_{\partial E} dh(\vec\nu)>0.$$
In particular $h$ is not the constant function and $h>0$ on $E$ by the maximum principle. We also have
$$\int_E |dh|^2\le \liminf_{R\to\infty} \int_E |dh_R|^2=\liminf_{R\to\infty} C(R)=C(\infty)<\infty.$$
We must show that $\inf_{E} h=0$. Let $v \in C^\infty_0(\overline{E})$ such that $v\ge 1$ on $\partial E$, we
let $v=h\varphi$ then using the Green formula we get :

\begin{equation*}
\begin{split}\int_E |dv|^2&=\int_E h^2|d \varphi|^2+2
\int_E h\varphi \langle dh, d\varphi\rangle+\int_E |dh|^2\varphi^2\\
&=\int_E h^2|d \varphi|^2+\frac12 \int_E \langle dh^2, d\varphi^2\rangle+\int_E |dh|^2\varphi^2\\
&=\int_E h^2|d \varphi|^2+\frac12 \int_E \Delta(h^2)\varphi^2-
 \int_{\partial E} \varphi^2 dh(\vec\nu)+\int_E |dh|^2\varphi^2\\
&\ge \int_E h^2|d \varphi|^2-\int_{\partial E}  dh(\vec\nu)\\
&\ge (\inf_{E} h)^2 C(\infty)+C(\infty)\\
\end{split}
\end{equation*}
Taking the infimum over such $v$ we obtain that 
$$0\ge (\inf_{E} h)^2 C(\infty)$$
hence $\inf_{E} h=0$

It remains to show that $iv)\Rightarrow iii)$. 
We re consider the preceding notation. Under the hypothesis $iv)$, we must show that $C(\infty)>0$.
Or by contraposition that if $C(\infty)=0$ then iv) can't be true. 
So we assume that $C(\infty)=0$, in that case we know that $h$ is the constant function $1$ and
we get for the functions $h_R$ :
$$\lim_{R\to\infty}\int_U h_R^2=\vol U$$ where as
$$\lim_{R\to\infty}\int_E |dh_R|^2=0.$$
That is iv) is not true

\endproof

\begin{prop}
If all ends of $M$ are non parabolic then 
$$\{0\} \rightarrow H^1_0(M)\rightarrow H^1_2(M). $$
In fact if $M$ has more than $k$ non parabolic ends then 
$$\dim\big\{h\in C^\infty(M), \Delta h=0\mbox{ and  }dh\in L^2\big\}\ge k.$$\end{prop}
The above inequality is due to P. Li and L-F Tam (\cite{LT}).
\proof
We remark that because $M$ contains at most one non parabolic end, the above proof of 
(\ref{para} {\it v)}$\Rightarrow$ {\it iv)}\, ) shows that
for any $U \subset M$ bounded open subset there is $C_U>0$ such that
$$\forall f\in  C^\infty_0(M), C_U\int_U f^2\le\int_M |df|^2.$$
We consider a closed form $\alpha\in C^\infty_0(T^*M)$ which is map to zero in $H^1_2(M)$.
Hence there is a sequence $v_k\in C^\infty_0(M)$ such that
$$\alpha=L^2-\lim_{k\to\infty } dv_k.$$
By the above remark, we know that there is $v\in C^\infty(M)$ such that  $dv=\alpha$ and
$$v=L^2_{loc}-\lim_{k\to\infty } v_k.$$
Let $E$ be a unbounded connected component of $M\setminus \supp \alpha$, $v$ is then
constant on $E$
$$v=v(E)\mbox{ on } E.$$ We want to show that this constant is zero. We can enlarge $M\setminus E$ and assume that
$E$ has smooth
boundary. On $E$ we have
$$\lim_{k\to\infty}\int_E|dv_k|^2=\int_E|\alpha|^2=0.$$
Choose $U$ a non empty bounded open subset of $E$, we know that for a certain constant $C_U>0$ we have
 the estimate
$$\int_E|dv_k|^2\ge C_U\int_U v_k^2,$$
When $k$ tends to $\infty$ we obtain
$$0\ge C_U\vol U\, v(E)^2.$$
Hence $v$ is zero on $E$ and $v$ has necessary compact support.
\endproof

\subsection{The two dimensional case}

In dimension $2$, a remarkable property of the space of $L^2$ harmonic $1$ forms is that it is an invariant
of conformal structure (or of the complex structure).
\subsubsection{Conformal invariance}
\begin{prop}
Let $\bar g=e^{2u}g$ be two conformally equivalent Riemannian metric
\footnote{The metric $\bar g$ and $g$ are not necessary complete.}
 on a smooth manifold $M^{2m}$
, then 
then we have the equality :
$$\cH^m(M,\bar g)=\cH^m(M,g).$$
\end{prop}
\proof When $\alpha\in \Lambda^kT^*_xM$, we have
$$|\alpha|^2_{\bar g}=e^{-2k u(x)} |\alpha|^2_{g}.$$
and $$d\vol_{\bar g}=e^{2m u}d\vol_{ g}.$$
As a consequence, the two Hilbert spaces 
$L^2(\Lambda^mT^*M, \bar g)$ and $L^2(\Lambda^mT^*M, g)$ are isometric, hence by definition the two space
$Z^m_2(M,\bar g)=\{\alpha\in L^2(\Lambda^mT^*M, \bar g), d\alpha=0\}$ and
$Z^m_2(M,\bar g)=\{\alpha\in L^2(\Lambda^mT^*M,  g), d\alpha=0\}$ are the same.
Moreover the orthogonal of $dC^\infty_0(\Lambda^{m-1} T^*M)$ is the same in $L^2(\Lambda^mT^*M,  g)$ or 
$L^2(\Lambda^mT^*M, \bar g)$. But by definition
$$\cH^m(M, \bar g)=Z^m_2(M,\bar g)\cap \left(dC^\infty_0(\Lambda^{m-1} T^*M)\right)^\perp.$$
is also equal to $\cH^m(M, g)$.

Another proof is to compute the codifferential $d^*_{\bar g}$ :
for $\alpha\in C^\infty(\Lambda^kT^*M)$ we obtain
$$d^*_{\bar g}\alpha=e^{-2u}\left(d^*_{ g}\alpha-2(m-k)\inte_{\grad u}\alpha\right).$$\endproof
\subsubsection{Application}
\begin{thm}\label{surface} Let $(S,g)$ be a complete connected Riemannian surface with finite topology (finite genus and finite
number of ends) then either 
\begin{itemize}
\item $\dim \cH^1(S,g)=\infty$
\item or $\dim \cH^1(S,g)<\infty$ and $M$ is conformally equivalent to a compact Riemannian surface $(\bar S,
\bar g)$ with a finite number of points removed :
$$(S,g)\simeq (\bar S\setminus \{p_1,...,p_k\}, \bar g)$$ and
$$\cH^1(S,g)\simeq \ima\left(H^1_0(S)\rightarrow H^1(S)\right)\simeq H^1(\bar S).$$
\end{itemize}
\end{thm}
\proof We know that a Riemannian surface with finite topology is necessary conformally 
equivalent to a compact Riemannian surface $(\bar S,\bar g)$ with a finite number of points and disks removed
:
$$(S,g)\simeq \left(\bar S\setminus \Big(\cup^{b}_{l=1} D_l\cup\{p_1,...,p_k\}\Big), \bar g\right)$$
Hence from our previous result, we have  
$$\cH^1(S,g)=\cH^1\left(\bar S\setminus \Big(\cup_{l=1}^b D_l\cup\{p_1,...,p_k\}\Big), \bar g\right).$$
We first show that if $b\ge 1$ then $\dim\cH^1(S)=\infty$.
Let $f\in C^\infty(\partial D_1)$, then on $\bar S\setminus D_1$ we can solve the Dirichlet problem :
$$\left\{\begin{array}{ll}
\Delta^{\bar g} u=0&\mbox{ on } \bar S\setminus D_1\\
u=f &\mbox{ on } \partial D_1\\
\end{array}\right.$$
Then $u\in C^\infty(\bar S\setminus D_1)$ hence $du\in L^2(T^*\left(\bar S\setminus D_1 \right),\bar g)$ is
closed and coclosed hence $du\in \cH^1(\bar S\setminus D_1,\bar g)$ its restriction to 
$S=\bar S\setminus \Big(\cup_{l=1}^{b} D_l\cup\{p_1,...,p_k\}\Big)$ is also $L^2$ closed and co closed, hence
we have build a linear map
$$f\in C^\infty(\partial D_ 1)\mapsto du\big|_{S}\in\cH^1(S,g)=\cH^1(S,\bar g)\ ,$$
the kernel of this map is the set of constant functions hence we have proved that
$$\dim\cH^1(S,g)=\infty.$$
Now we assume that $b=0$ that is $$(S,g)\simeq (\bar S\setminus \{p_1,...,p_k\}, \bar g)\ ,$$ 
we have
$$\cH^1(S,g)=\cH^1(\bar S\setminus \{p_1,...,p_k\}, \bar g).$$
The main point is that when $\alpha\in \cH^1(\bar S\setminus \{p_1,...,p_k\}, \bar g)$ then
$\alpha$ extends across the point $\{p_1,...,p_k\}$.  
This is a direct consequence of the following lemma
\begin{lem}
Let $\bD$ be the unit disk, and let 
$\alpha\in L^2\big(T^*(\bD\setminus\{0\})\big)$ satisfying the equation 
$(d+d^*)\alpha=0$, then 
$\alpha$ extends smoothly across $0$ that is 
$$\cH^1(\bD\setminus\{0\})=\cH^1(\bD).$$\end{lem}\endproof
In this lemma the metric on the disk can be considered as the Euclidean one, as a matter of fact any other
Riemannian metric on $\bD$ is conformally equivalent to the flat metric. We offer two proof of this result.

{\it First proof of the lemma.}  Let $\alpha\in\cH^1(\bD\setminus\{0\})$, we are going to prove that 
$ (d+d^*)\alpha=0$ holds weakly on $\bD$, 
then because the operator $(d+d^*)$ is elliptic it will hold strongly by elliptic regularity.
Consider $\varphi\in C^\infty_0(\bD)\oplus C^\infty_0(\Lambda^2T^*\bD) $
we must show that 
$$\langle \alpha, (d+d^*)\varphi\rangle=0.$$
We consider the following sequence of cutoff functions
$$\chi_n(r,\theta)=\left\{\begin{array}{ll}
0& \mbox{ when }  r\le 1/n^2\\
-\frac{\log(r/n^2)}{\log n}&\mbox{ when } 1/n^2\le r\le 1/n\\
1&\mbox{ when } r\ge 1/n\\
\end{array}\right.$$
By hypothesis 
$$\langle \alpha, (d+d^*)(\chi_n\varphi)\rangle=0.$$
But $(d+d^*)(\chi_n\varphi)=\chi_n (d+d^*)\varphi+d\chi_n\wedge \varphi+\inte_{\grad \chi_n} \varphi$.
We have then 
$$0=\langle \alpha, (d+d^*)(\chi_n\varphi)\rangle=
\langle \alpha, \chi_n(d+d^*)\varphi)\rangle+
\langle \alpha, d\chi_n\wedge \varphi+\inte_{\grad \chi_n} \varphi\rangle$$

But when $n\to\infty$, the first term in the right hand side goes to $\langle \alpha, (d+d^*)\varphi\rangle$ :
$$\lim_{n\to\infty} \langle \alpha, \chi_n(d+d^*)\varphi)\rangle=\langle \alpha, (d+d^*)\varphi)\rangle$$
and the second term is estimate as follow :
$$\Big|\langle \alpha, d\chi_n\wedge \varphi+\inte_{\grad \chi_n} \varphi\rangle\Big|\le
 \|\alpha\|_{L^2}\,\|\varphi\|_{L^\infty }\|d\chi_n\|_{L^2}$$
 But a direct computation shows that 
 $$\|d\chi_n\|_{L^2}^2=\frac{2\pi}{\log n}.$$ Hence the result. \hfill$\square$
 
{\it Second proof of the lemma.} Let $\alpha\in \cH^1(\bD\setminus\{0\})$.
We first show that $\alpha$ is exact. That is we'll prove that
$$\int_{\bS^1}\alpha=0$$
Let $c=\int_{\bS^1}\alpha$, this integral doesn't depend on the radius ; hence we get by Cauchy-Schwarz
inequality we get :
$$c^2\le 2\pi r^2 \int_0^{2\pi} |\alpha|^2(r,\theta) d\theta $$
Dividing by $r$ and integrating over $[\varepsilon,1]$ we get :
$$-c^2\log \varepsilon\le 2\pi \int_{\bD} |\alpha|^2.$$
Hence letting $\varepsilon\to 0+$ we obtain $c=0$. So there is an harmonic function $f$ on
$\bD\setminus\{0\}$ 
 such that
$$\alpha=df.$$ 
We have 
$$f(r,\theta)=\sum_{k\in \Z} u_k(r) e^{ik\theta}$$
where $u_k(r)=A_k r^{|k|}+B_k r^{|k|}$ when
$k\not=0$ and
$u_0(r)=A_0+B_0\log(r)$, moreover
$$\int_{\bD} |df|^2=\sum_{k\in \Z} \int_0^1\left[ |u'_k|^2+k^2|u_k|^2\right] rdr<\infty.$$
This implies that $B_k=0$ for all $k\in \Z$· 
hence $f$ is smooth at zero.\hfill$\square$

\subsection{Bibliographical hints}

For the second lecture, we warmly recommend the reading of
the paper of A. Ancona \cite{Ancona} about non parabolicity, we also recommend the  beautiful survey 
by A. Grigor'yan \cite{Gri} on non parabolicity, stochastic completeness and harmonic functions. The paper of P.
Li-F.Tam \cite{LT} about the space of harmonic functions with $L^2$ differential has also to be read.

\section{$L^2$ cohomology of conformally compact manifold}
\subsection{The geometric setting}
Let $\overline{M}$ be a compact smooth manifold with boundary $N=\partial\overline{M}$. On 
$M=\mbox{int}(\overline{M})=\overline{M}\setminus N$, we say that a Riemannian metric $g$ is 
{\it conformally compact} if 
$$g=\frac{\bar g}{y^2}$$
where $\bar g$ is a smooth Riemannian metric on $\bar M$ (hence non complete) and 
$y\,:\, \overline{M}\rightarrow \R_+$ is  smooth defining function for $N$ that is 
$$\left\{
\begin{array}{ll}y^{-1}\{0\}=N& \\
dy\not=0& \mbox{ along } N\\
\end{array}\right.$$
Such a metric is complete. The most famous example is the hyperbolic metric in the ball model :
$$g_{\rm hyp}=\frac{4 \|dx\|^2}{(1-\|x\|^2)^2}$$ where
$\|dx\|^2$ is the Euclidean metric on the unit ball of $\R^n$ and 
$$y=\frac{1-\|x\|^2}{2}$$ is a smooth defining function for $\bS^{n-1}=\partial\bB^n$.

Conformally compact Riemannian metrics have a great interest in theoretical physics with regards to the
AdS/CFT correspondence \cite{Bi}. In fact these class of metric have been first study by C. Fefferman and R. Graham
(\cite{GF}).

\subsection{The case where $k=\dim M/2$.}
\begin{lem}
Assume that $(M,g)$ is a conformally compact Riemannian manifold of dimension $2k$ then
$$\dim H^k_2(M)=\infty.$$
\end{lem}
\proof Because $g$ and $\bar g$ are conformally equivalent we know that 
$$H^k_2(M)=\cH^k(M,g)=\cH^k(M,\bar g).$$
 Hence this lemma will a consequence of the following result
 \begin{lem}
 Assume that $(\overline{M}, \bar g)$ is compact Riemannian manifold with smooth boundary, then for every
 $k\not=0, \dim M$, we have 
 $$\dim\{\alpha\in C^\infty(\Lambda^kT^*\overline{M}), d\alpha=d^*\alpha=0\}]=\infty.$$\end{lem}
 When $k=0$, we know that $\cH^0(\overline{M})$ consist of locally constant function hence is finite
 dimensional. Moreover, this equality for $k=1$ can be proved by the same argument of the proof of
  the first assertion of (\ref{surface}), when $\overline{M}$ is connected, we have a 
  linear map from the space of smooth function on the boundary of $\overline{M}$ to $\cH^1(\overline{M})$
  which associated to each $f\in C^\infty(\partial\overline{M})$ the differential of its harmonic extension
  to $\overline{M}$. The kernel of this map have dimension $1$ hence when $n\ge 2$ we obtain that 
 $\dim \cH^1(\overline{M})=\infty.$
 
 A general proof using pseudo differential calculus shows that the solution space of an elliptic operator
 on a compact manifold with boundary has infinite dimension.  The proof given below is more elementary, we
 will only used the {\it unique continuation principle}:
 \begin{prop} If $N$ is a connected Riemannian manifold and $U\subset N$ is a non empty open subset of $N$.
 Then if $\alpha,\beta\in \cH^k(N)$ satisfy $\alpha=\beta $ on $U$ then
 $$\alpha=\beta\mbox{  on  } N.$$
 \end{prop}
 We consider $D=\overline{M}\#_{\partial \overline{M}}\overline{M} $ the double of $\overline{M}$, 
 and we
consider
$$N_l=D\# \Big(\#_l \bT^n\Big)$$ the manifold obtained by making the connected sum of $D$ with $l$ copy of a
$n$-torus $\bT^n$ on the second copy of $\overline{M}\subset D$. 
We endowed $N_l$ with a smooth Riemannian metric which coincide with $\bar g$ on the first copy. 
\begin{center}\fbox{\includegraphics[height=3cm]{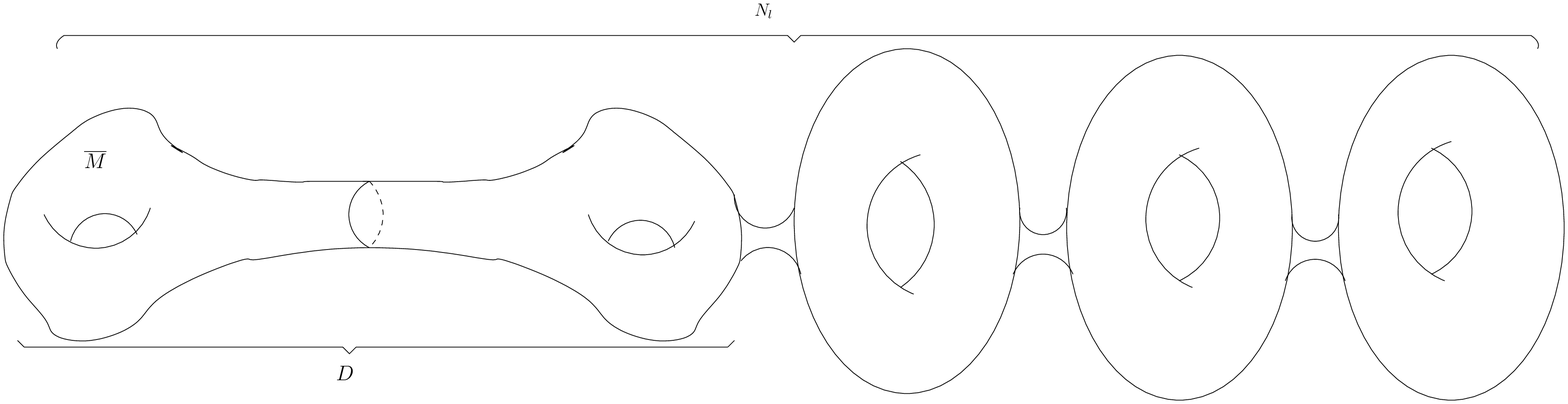}  }
\end{center}
Then by the unique continuation principle we know that the restriction map $\cH^k(N_l)\rightarrow \cH^k(\bar
M)$ on the first copy  is injective
hence
$$\dim \cH^k(\bar M,\bar g)\ge  b_k(N_l)\,.$$
However when $k\not= 0, \dim \bar M$ we have
$b_k(N_l)=b_k(D)+l b_k(\bT^n)$, hence
letting $l\to\infty$ we obtain the desired result.
\endproof

\subsection{A reduction to exact metric}
The spaces of reduced $L^2$ cohomology depends only on the $L^2$ structures, hence if $g_1$ and $g_2$
are two Riemannian metric on a manifold $M$ ($g_1,g_2$ need not to be complete) which are quasi isometric
that is  for a certain constant $C$
$$\frac{g_1}{C} \le g_2 \le C g_1\ , $$ then clearly the Hilbert spaces 
$L^2(\Lambda^k T^*M,g_1) $ and $L^2(\Lambda^k T^*M,g_2) $ are the same with equivalent norms. Hence the
quotient space defining reduced $L^2$ cohomology are the same\footnote{Indeed the space $Z^k_2(M)$
depends only of the topology of $L^2(\Lambda^k T^*M,g)$ : on 
$L^2(\Lambda^k T^*M,g)\times C^\infty_0(\Lambda^{k+1}T^*M)$ the bilinear form $(\alpha,\beta)\mapsto
\langle\alpha,d^*\beta\rangle$ does not depends on $g$.} , that is
$$H^k_2(M,g_1)=H^k_2(M,g_2).$$ In particular we have the following :
 \begin{prop}
 If $(M,g)$ is a conformally compact Riemannian metric, then $H^k_2(M,g)$ doesn't depend
 on the conformally compact metric on $M$.
 \end{prop}
 \proof As a matter of fact, if
 $g=\frac{\bar g}{y^2}$ and 
 $g_1=\frac{\bar g_1}{y_1^2}$ are two conformally compact Riemannian metric on $M$, then because 
 $\bar g$ and $\bar g_1$ are smooth Riemannian metric on a compact manifold, we know that there is a constant
 $C_1$ such that 
 $$\frac{\bar g_1}{C_1} \le\bar g \le C_1\bar g_1.$$
 Moreover there is always a smooth function $u\,:\, \overline{M}\rightarrow \R$ such that 
 $$y=e^u y_1$$
 hence there is a constant $C_2$ such that
 $$\frac{y_1}{C_2} \le y \le C_2 y_1.$$
 Eventually we obtain for $C=C_1 C_2^2$, $$C^{-1} g_1 \le g \le C g_1.$$
 \endproof
 
 An simpler example of conformally compact Riemannian metric are the so-called exact conformally compact
 metric : let $h$ be a smooth Riemannian metric on $N=\partial\overline{M}$ and 
 $y\,:\, \overline {M}\rightarrow\R_+$ a boundary defining function, if $\varepsilon>0$ is sufficiently small
 then 
 $$\frac{|dy|^2+h}{y^2}=dr^2+e^{2r} h, r=-\log y$$
 is a conformally compact metric on the collar neighborhood 
 $\{y< \varepsilon\}$ of $N\subset \overline{M}$. Any conformally compact Riemannian metric having such an
 expression near $N$ will be called exact.
\subsection{A Rellich type identity and its applications}
\subsubsection{The formula}
Let $\Omega$ be a Riemannian manifold with smooth boundary $\partial \Omega$ 
\begin{lem} \label{formula} Let $X$ be a vector field on $\Omega$ and $\alpha\in C^\infty(\Lambda^k T^*\Omega)$ then
\begin{equation*}
\begin{split}
\cL_X\alpha=\left.\frac{d}{dt}\right|_{t=0}\left(\Phi_X^t\right)^*\alpha&=d(\inte_X\alpha)+\inte_X(d\alpha)\\
&=\nabla X.\alpha+\nabla_X\alpha\ ,\\
\end{split}
\end{equation*}
where $\Phi_X^t$ is the flow associated to $X$ and 
for $(E_1,...,E_n)$ a local orthonormal frame and $(\theta^1,...,\theta^n)$ its dual frame then 
$$\nabla X.\alpha=\sum_{i,j}\theta^i\wedge \theta^j\,(\nabla_{E_i}X,E_j) \inte_{E_i}\inte_{E_j} \alpha
=\sum_{i=1}^n\theta^i\wedge(\inte_{\nabla_{E_i}X}\alpha).$$
\end{lem}
\proof We used the formula (\ref{expd})
$$d=\sum_{i=1}^n\theta^i\wedge\nabla_{E_i} \ .$$
Hence
\begin{equation*}
\begin{split}
d(\inte_X\alpha)&=\sum_{i=1}^n\theta^i\wedge\nabla_{E_i}(\inte_X\alpha)\\
&=\sum_{i=1}^n\theta^i\wedge(\inte_{\nabla_{E_i}X}\alpha)+
\sum_{i=1}^n\theta^i\wedge(\inte_{X}\nabla_{E_i}\alpha)\\
&=\nabla X.\alpha+\sum_{i=1}^n\theta^i\wedge(\inte_{X}\nabla_{E_i}\alpha)\\\\
\end{split}
\end{equation*}
and
\begin{equation*}
\begin{split}
\inte_X(d\alpha)&=\sum_{i=1}^n\inte_X(\theta^i\wedge\nabla_{E_i}\alpha)\\
&=\sum_{i=1}^n\inte_X(\theta^i)\,\nabla_{E_i}\alpha
-\sum_{i=1}^n\theta^i\wedge\inte_X(\nabla_{E_i}\alpha)\\
&=\nabla_X\alpha-\sum_{i=1}^n\theta^i\wedge\inte_X(\nabla_{E_i}\alpha)\\
\end{split}
\end{equation*}
Hence the result.\endproof

\begin{cor}Let $\vec\nu\,:\, \partial\Omega\rightarrow T\Omega$ be the inward unit vector field, then for
$\alpha\in C^\infty_0(\Lambda^kT^*\overline{\Omega})$ :
\begin{equation*}
\begin{split}
\int_\Omega (\nabla X.\alpha,\alpha)+\frac12 \dive X\, |\alpha|^2=&\int_\Omega
(\inte_X\alpha,d^*\alpha)+(\inte_Xd\alpha,\alpha)\\
&\ \ +\int_{\partial\Omega}\left[
\frac12 ( X,\vec\nu) |\alpha|^2-(\inte_{\vec\nu}\alpha,\inte_{X}\alpha)\right] d\sigma.\\
\end{split}
\end{equation*}
\end{cor}
\proof The lemma (\ref{formula}) implies that 
$$\int_\Omega \left[(\nabla X.\alpha,\alpha)+(\nabla_X\alpha,\alpha)\right]d\vol_g=
\int_\Omega \left[(d(\inte_X\alpha),\alpha)+(\inte_X(d\alpha),\alpha) \right]d\vol_g.$$
Then the equality follows directly from the following two Green's type formulas :
$$\int_\Omega(\nabla_X\alpha,\alpha)d\vol_g=\frac12\int_\Omega
X.(\alpha,\alpha)d\vol_g= 
\frac12\int_\Omega \dive X |\alpha|^2\,d\vol_g-\frac12\int_{\partial\Omega} ( X,\vec\nu) |\alpha|^2d\sigma$$
and 
$$ \int_\Omega (d(\inte_X\alpha),\alpha)d\vol_g=\int_\Omega (\inte_X\alpha,d^*\alpha)d\vol_g-
\int_{\partial\Omega}(\inte_{\vec\nu}\alpha,\inte_{X}\alpha) d\sigma.$$
\endproof
This integration by part formula is due to Donnelly-Xavier in \cite{DX}, this formula has been used and
refined
by many authors (see for instance \cite{EF},\cite{Kasue}).
\subsection{Application to conformally compact Riemannian manifold}
We apply here this formula to $\Omega=]0,\infty[\times N$ endowed with the metric
$dr^2+e^{2r} h$ where $h$ is a smooth Riemannian metric on the compact manifold $N$. We choose the vector
field $$X=-\frac{\partial}{\partial r}.$$ The curves $r\mapsto (r,\theta)\in \Omega$ are geodesic hence
$$\nabla_X \frac{\partial}{\partial r}=0 ,$$ moreover it is not hard to verify that for $v\in TN$
$$ \nabla_vX=-v\ \footnote{this comes from the fact that that level set of $r$ are totally umbilical.}.$$
A direct computation shows that for
$\alpha\in \Lambda^kT^*\Omega$
$$(\nabla_X\alpha,\alpha)=-k |\alpha|^2+|\inte_X\alpha|^2\mbox{ and } \dive X=(n-1).$$
Hence 
$$(\nabla_X\alpha,\alpha)+\frac12 \dive X
|\alpha|^2=\left(\frac{n-1}{2}-k\right)|\alpha|^2+|\inte_X\alpha|^2$$
and we obtain that for all $\alpha\in C^\infty_0(\Lambda^kT^*\overline{\Omega})$ :
\begin{equation}\label{estim1}\begin{split}
\int_\Omega\left[(\inte_X\alpha,d^*\alpha)+(\inte_Xd\alpha,\alpha) \right]d\vol_g\ge
\left(\frac{n-1}{2}-k\right)&\int_\Omega |\alpha|^2d\vol_g\\
&+\int_{\partial\Omega}
\left[\frac12 |\alpha|^2-|\inte_{X}\alpha|^2\right] d\sigma.\\
\end{split}\end{equation}
\begin{cor}\label{vanishabs} Assume that $k\le \frac{n-1}{2}$ and that
$\alpha\in L^2(\Lambda^kT^* \Omega)$ satisfies 
$$d\alpha=d^*\alpha=0$$
and 
$$\inte_ {\vec  \nu} \alpha=0\mbox{ along }  \partial \Omega $$
then 
$\alpha=0$. Hence by (\ref{hodgeabs}) we have that 
$H_2^k(\Omega)=\{0\}$ for  $k<n/2$.
\end{cor}
\proof Let $\rho\in C^\infty_0(\R_+)$ having support in $[0,1]$ and such that $\rho=1$ near $0$ and let
$\rho_N(t)=\rho(t/N)$. We apply the
inequality (\ref{estim1}) to $\alpha_N=\rho_N\alpha$ :
$$\int_\Omega(\inte_Xd\alpha_N,\alpha_N) d\vol_g\ge
\left(\frac{n-1}{2}-k\right)\int_\Omega |\alpha_N|^2d\vol_g+
\frac12\int_{\partial\Omega} |\alpha|^2 d\sigma.$$
But 
$$\int_\Omega(\inte_Xd\alpha_N,\alpha_N) d\vol_g=\int_\Omega \frac{1}{N} \rho'\left(\frac{r}{N}\right)
|\alpha|^2 d\vol_g\le \frac{1}{N}\|\rho'\|_{L^\infty}\int_\Omega|\alpha|^2 d\vol_g.$$
Hence letting $N$ going to infinity, we obtain that
$$0=\left(\frac{n-1}{2}-k\right)\int_\Omega |\alpha|^2d\vol_g+
\frac12\int_{\partial\Omega} |\alpha|^2 d\sigma.$$
Hence if $k<(n-1)/2$, we obtain $\alpha=0$. When $k=(n-1)/2$, we obtain that
$\alpha=0$ along $\partial \Omega$. Then the unique continuation property for solution of elliptic 
operator of order $1$ implies that $\alpha=0$.\endproof
\begin{cor}  \label{gap}
Assume that $k<(n-1)/2$, then $\forall \alpha\in C^\infty_0(\Lambda^kT^*\Omega)$,

$$\big\|d\alpha\big\|_{L^2}^2+\big\|d^*\alpha\big\|_{L^2}^2\ge
 \frac12\left(\frac{n-1}{2}-k\right)^2\big\|\alpha\big\|_{L^2}^2 .$$
 \end{cor}
 \proof We apply the estimate (\ref{estim1}) to $\alpha\in C^\infty_0(\Lambda^kT^*\Omega)$, then there is no
 boundary term and we get the inequality 
 $$\int_\Omega\left[(\inte_X\alpha,d^*\alpha)+(\inte_Xd\alpha,\alpha) \right]d\vol_g\ge
\left(\frac{n-1}{2}-k\right)\int_\Omega |\alpha|^2d\vol_g .$$
But with the Cauchy-Schwarz inequality, we obtain :
\begin{equation*}
\begin{split}
\int_\Omega\left[(\inte_X\alpha,d^*\alpha)+(\inte_Xd\alpha,\alpha) \right]d\vol_g
&\le
 \big\|d^*\alpha\big\|_{L^2}\,\big\|\alpha\big\|_{L^2}
 +\big\|d\alpha\big\|_{L^2}\,\big\|\alpha\big\|_{L^2}\\
 &\le
  \sqrt{2}\left[\big\|d\alpha\big\|_{L^2}^2+\big\|d^*\alpha\big\|_{L^2}^2 \right]^{1/2}
  \big\|\alpha\big\|_{L^2}.\\
 \end{split} \end{equation*}
  Hence the result.
 \endproof
 \begin{rem}\label{vanishrel}
 Using the vector field $-X$, 
 the reader can check that the estimate of the corollary (\ref{gap})  is also true for $k>(n+1)/2$. Also with
 (\ref{hodgerel}), the reader can also prove that $H_2^k(\Omega,\partial \Omega)=\{0\}$ for $k>n/2$.
 \end{rem}
\subsection{The spectrum of the Hodge-deRham Laplacian}
\subsubsection{The essential spectrum} Let $A\,:\, \cD(A)\rightarrow H$ be a selfadjoint
operator on a Hilbert space, then the spectrum of $A$ is the subset $\Spec A$ of $\C$ consisting of those
$z\in \C$ such that $A-z\Id$ has not a bounded inverse. Because $A$ is self adjoint, we have
$\Spec A\subset \R$. Moreover if $\cC\subset \cD(A)$ is a core for $A$ then we have that 
$\lambda\in \R$ belongs to the spectrum of $A$, if and only if there is a sequence $(\varphi_n)_n$ of element
of $\cC$ such that 
$$\left\{\begin{array}{l}\|\varphi_n\|=1\\
\lim_{n\to\infty}\|A\varphi_n-\lambda\varphi_n\|=0.\\
\end{array}\right.$$
 We can separate the spectrum of $A$ in two parts the discrete part and the essential part :
 $$\Spec A= \Spec_d A\cup\Spec_e A,$$
 where $\Spec_d A$ is the set of isolated point in $\Spec A$ which are eigenvalue with finite multiplicity
 and $\Spec_e A$ is the set of non-isolated point in $\Spec A$ or of eigenvalue with infinite multiplicity.
 We have the following characterization of the essential spectrum :
  a real number $\lambda$ belongs to the essential spectrum of $A$ if and only if

 there is a sequence $(\varphi_n)_n$ in $\cD(A)$ with 
$$\left\{\begin{array}{l}\|\varphi_n\|=1\\
\lim_{n\to\infty}\varphi_n=0\mbox{ weakly in }H\\
\lim_{n\to\infty}\|A\varphi_n-\lambda\varphi_n\|=0.\\
\end{array}\right. \\
$$
Another useful characterisation of the complementary of the essential spectrum is the following :
\begin{prop}\label{green}

$\lambda\not\in \Spec_e A$ if and only if there is a bounded operator
$G\,:\, H\rightarrow \cD(A)$ such that $(A-\lambda\Id)G-\Id$ and $G(A-\lambda\Id)-\Id$ are compact operators.
Moreover the operator $G$ can be chosen so that $(A-\lambda\Id)G-\Id=G(A-\lambda\Id)-\Id$ is the orthogonal
projection on $\ker(A-\lambda\Id)$.
\end{prop}
From these properties, it is not hard to verify that the essential spectrum is stable by compact perturbation :
\begin{thm} Let $A,B\,:\, \cD(A)=\cD(B)\rightarrow H$ be two self-adjoint operator such that 
$(A+i)^{-1}-(B+i)^{-1}$ is a compact operator then $A$ and $B$ have the same essential spectrum : $$\Spec_e(A)=\Spec_e(B).$$
\end{thm}

\subsubsection{Case of  the Hodge-deRham Laplacian of complete Riemannian manifold}
Let $(M,g)$ be a complete Riemannian manifold, then the operator $\Delta=dd^*+d^*d\,:\,
C^\infty_0(\Lambda^kT^*M)\rightarrow L^2(\Lambda^kT^*M)$ has a unique self-adjoint extension to  
$L^2(\Lambda^kT^*M)$ which also denoted by $\Delta$ with domain
$$\cD(\Delta)=\{\alpha\in L^2(\Lambda^kT^*M), \Delta\alpha\in L^2\}$$
That is 
$\alpha\in \cD(\Delta)$ if and only if there is a constant $C>0$ such that 
$$\forall \varphi\in C^\infty_0(\Lambda^kT^*M),
\ |\langle\alpha,\Delta\varphi\rangle|\le C \|\varphi\|_{L^2}.$$
In fact we can prove that 
$$\cD(\Delta)=\{\alpha\in L^2(\Lambda^kT^*M), d\alpha\in L^2,d^*\alpha\in L^2, dd^*\alpha\in L^2,
d^*d\alpha\in L^2\}.$$
Moreover, $C^\infty_0(\Lambda^kT^*M)$ is dense in $\cD(\Delta)$ for the graph norm :
$\alpha\mapsto \|\alpha\|_{L^2}+\|\Delta\alpha\|_{L^2}$. We have the following very useful result of H.
Donnelly and I.Glazman 
(see \cite{Ang,Ba,Do,Glaz})
\begin{thm}\label{gla} Zero is not in the essential spectrum of $\Delta$ $(0\not\in \Spec_e(\Delta))$ if and only if
there is a compact set $K\subset M$ and a constant $\varepsilon>0$ such that
$$\forall \varphi\in C^\infty_0(\Lambda^kT^*(M\setminus K)),
 \|d\alpha\|_{L^2}^2+\|d^*\alpha\|_{L^2}^2=\langle\alpha,\Delta\alpha\rangle\ge 
 \varepsilon \|\alpha\|^2_{L^2}.$$
\end{thm}
\begin{rem}
When $\Omega \subset M$ is a open set with smooth compact boundary $\partial \Omega$, we can introduce two
operators :
\begin{itemize}
\item The Laplacian with the relative boundary condition :
$$\cC_{rel}=\{\alpha\in C^\infty_0(\Lambda^kT^*\overline{\Omega}),\mbox{ such that
}\iota^*\alpha=\iota^*(d^*\alpha)=0\}$$
then $\Delta\,:\,\cC_{rel}\rightarrow L^2$ has a unique selfadjoint extension $\Delta_{rel}$.
\item The Laplacian with the absolute boundary condition :
$$\cC_{abs}=\{\alpha\in C^\infty_0(\Lambda^kT^*\overline{\Omega}),\mbox{ such that
}\inte_{\vec\nu}\alpha=\inte_{\vec\nu}d\alpha=0\}$$
then $\Delta\,:\,\cC_{abs}\rightarrow L^2$ has a unique selfadjoint extension $\Delta_{abs}$.
\end{itemize}
Moreover when $M\setminus \Omega$ is a compact set  we have
$$0\not\in \Spec_e \Delta\Leftrightarrow 0\not\in \Spec_e \Delta_{abs}\Leftrightarrow 
0\not\in \Spec_e \Delta_{rel}.$$
\end{rem}
\subsection{Applications to conformally compact manifolds}
\subsubsection{The essential spectrum}
\begin{prop} Let $(M,g)$ be a conformally compact Riemannian manifold endowed with an exact metric, then for 
$k\not\in\left\{\frac{n-1}{2},\frac{n}{2},\frac{n+1}{2}\right\}$, zero is not in the essential spectrum
of the Hodge deRham Laplacian acting on $k-$forms.
\end{prop}
\proof This is a direct corollary of the result (\ref{gla}) of I.Glatzmann and H.Donnelly and of the corollary (\ref{gap}).
\endproof
\begin{cor}\label{primitive}
Assume again that $(M,g)$ be a conformally compact Riemannian manifold endowed with an exact metric. For
$k\le (n-1)/2$, we consider $\alpha\in Z^k_2(M)$ which is zero in $H^k_2(M)$ then there is $\beta\in
L^2(\Lambda^{k-1}T^*M)$
such that $$\alpha=d\beta.$$
\end{cor}
\proof By hypothesis, we know that there is a sequence $\varphi_l \in C^\infty_0(\Lambda^{k-1}T^*M)$ such
that
$$\alpha=L^2-\lim_{l\to\infty} d\varphi_l.$$
Because $k-1<(n-1)/2$, from (\ref{green}) we know that there is a bounded operator
$$G\,:\,L^2(\Lambda^{k-1}T^*M)\rightarrow L^2(\Lambda^{k-1}T^*M) $$ such that 
$\forall \varphi\in C^\infty_0(\Lambda^{k-1}T^*M)$
$$\Delta G\varphi=\varphi-h(\varphi),$$
where $h(\varphi)$ is the orthogonal projection on the kernel of $\Delta$ that is on $\cH^k(M)$.
Hence for $\psi_l=G\varphi_l$ we obtain
$$\varphi_l=h(\varphi_l)+dd^*\psi_l+d^*d\psi_l$$
and 
$$\alpha=L^2-\lim_{l\to\infty} dd^*d\psi_l.$$
We let $\eta_l=d^*d\psi_l$, we have
$d^*\eta_l=0$ and $d\eta_l=d\varphi_l$, hence, by elliptic regularity,
$\eta_l$ is smooth. In particular, $d\eta_l\in L^2$ and $d^*\eta_l\in L^2$ and
$h(\eta_l)=0$.
We have 
\begin{equation*}
\begin{split}
\|\eta_l-\eta_k\|_{L^2}^2&=\langle \Delta G(\eta_l-\eta_k),(\eta_l-\eta_k)\rangle\\
&=\langle dG (\eta_l-\eta_k),d(\eta_l-\eta_k)\rangle
                  + \langle d^* G(\eta_l-\eta_k),d^*(\eta_l-\eta_k)\rangle\\
&\le  \|dG (\eta_l-\eta_k)\|_{L^2}\|d(\eta_l-\eta_k)\|_{L^2}\\
&\le \left[\, \langle \Delta G (\eta_l-\eta_k),G(\eta_l-\eta_k)\rangle\,\right]^{1/2}\|d(\eta_l-\eta_k)\|_{L^2}\\
&\le \left(\|\eta_l-\eta_k\|_{L^2}\|G(\eta_l-\eta_k)\|_{L^2} \right)^{1/2}\|d(\eta_l-\eta_k)\|_{L^2}\\
&\le C\|(\eta_l-\eta_k)\|_{L^2}\|d(\eta_l-\eta_k)\|_{L^2}.\\
\end{split}
\end{equation*}
where $C^2$ is the operator norm of $G$.
Hence 
$$\|\eta_l-\eta_k\|_{L^2}\le C\|d(\eta_l-\eta_k)\|_{L^2}=C\|d(\varphi_l-\varphi_k)\|_{L^2}.$$
The sequence $(d\varphi_l)_l$ is a Cauchy sequence in $L^2$ hence $(\eta_k)_k$ is also a Cauchy sequence in
$L^2$ converging to some $\beta\in L^2$ and
we have
$$\alpha=d\beta.$$
\endproof
\begin{rem}
This proof also shows that the primitive $\beta$ obtained satisfies the equation $d^*\beta=0$. Hence if 
$\alpha$ is smooth then $\beta$ will be also smooth.
\end{rem}
\subsubsection{At infinity}
With the same method, and because we have the vanishing result for the $L^2$ cohomology (\ref{vanishabs}),
 we have
\begin{prop}\label{primitive2} Let $\Omega=]0,\infty[\times N$ endowed with the warped product metric
$(dr)^2+e^{2r} h$, then for $k\le (n-1)/2$ and 
$\alpha\in Z^k_2(\Omega)$ there is $\beta\in L^2(\Lambda^{k-1}T^*\Omega)$ such that
$$\alpha=d\beta.$$
\end{prop}
\proof We offer a proof which is more elementary than the one passing through the 
essential spectrum and the vanishing result, this proof has an independent interest; the argument comes from
an article by P. Pansu (\cite{Pan}) where such techniques were used in order to obtain clever negative pinching results.
We introduce the vector field $T=-X=\frac{\partial}{\partial r}$ and we consider its flow
$$\Phi^t(r,\theta)=(r+t,\theta).$$
When $\alpha=dr\wedge\alpha_1+\alpha_2\in C^\infty(\Lambda^kT^*\Omega)$ where 
$\alpha_1\in C^\infty(\R_+,C^\infty(\Lambda^{k-1}T^*N))$ and $\alpha_2\in
C^\infty(\R_+,C^\infty(\Lambda^{k}T^*N))$, we easily obtain the estimate :
\begin{equation}\label{estimee}
\begin{split}
\left|\left(\Phi^t\right)^*\alpha\right|^2(r,\theta)&=e^{2(k-1)t}\left|\alpha_1\right|^2(r+t,\theta)+
e^{2kt}\left|\alpha_2\right|^2(r+t,\theta)\\
&\le e^{2kt}\left|\alpha\right|^2(r+t,\theta)\\
\end{split}
 \end{equation}

But when $\alpha\in Z^k_2(\Omega)$, the Cartan formula says that 
\begin{equation}\label{Cartan}\left(\Phi^t\right)^*\alpha-\alpha=d\beta_t. \end{equation}
where 
$$\beta_t=\int_0^t \left(\Phi^s\right)^*\left(\inte_T\alpha\right) ds.$$
From our estimate (\ref{estimee}), when $k\le (n-1)/2$ we easily obtain
\begin{equation}\label{decay1}
\begin{split}
\left\|\left(\Phi^t\right)^*\alpha\right\|_{L^2}^2&
\le \int_{]0,\infty[\times N} e^{2kt}\left|\alpha\right|^2(r+t,\theta) e^{(n-1)r} drd\vol_h(\theta)\\
&\le
\int_{]0,\infty[\times N} \left|\alpha\right|^2(r+t,\theta) e^{(n-1)(r+t)} drd\vol_h(\theta)\\
&\le\int_{]t,\infty[\times N} \left|\alpha\right|^2(r,\theta) e^{(n-1)r}
drd\vol_h(\theta)=\|\alpha\|^2_{L^2(\Omega_t)}\ ,\\
\end{split}
 \end{equation}
 where $\Omega_t=]t,\infty[\times N$. Hence
 $$L^2-\lim_{t\to+\infty} \left(\Phi^t\right)^*\alpha=0.$$
 
 Moreover for $t\ge 0$, we have 
\begin{equation}\label{decay2}
\begin{split} 
|\beta_t|(r,\theta)\le\int_0^t |\alpha|(r+s,\theta)e^{(k-1)s} ds\\
&\le e^{-(k-1)r}\int_r^\infty |\alpha|(s,\theta)e^{(k-1)s} ds.\\
\end{split}
 \end{equation}
 We need the following lemma
 \begin{lem}\label{continuite} Let $v\in C^\infty_0(\R_+)$ and let
 $$Mv(r)= e^{-(k-1)r}\int_r^\infty v(s)e^{(k-1)s} ds$$
 then for $k-1<(n-1)/2$ we have
 $$\left(\frac{n-1}{2}-(k-1)\right)^2\int_0^\infty |Mv|^2(t) e^{(n-1)t} dt\le  \int_0^\infty |v|^2(t) e^{(n-1)t} dt.$$
 That is the operator $M$ extends as a bounded operator in $L^2(\R_+, e^{(n-1)t} dt)$. 
 \end{lem}
 {\it Proof of the lemma (\ref{continuite}) . }
 Let $w(r)=e^{(k-1)r}Mv(r)$, then $w$ has compact support and
 $w'(r)=-e^{(k-1)r}v(r)$ and let $$w(r)=f(r) e^{-\left(\frac{n-1}{2}-(k-1)\right)r}$$ and 
 $\epsilon_k=n-1-2(k-1)>0$.
Then
 \begin{equation*}
\begin{split}
\int_0^\infty |v|^2(r) e^{(n-1)r} dr&=\int_0^\infty |w'|^2(r) e^{\epsilon_kr} dr\\
&=\int_0^\infty
 \left[ |f'|^2(r)-\epsilon_kf'f+\left(\frac{\epsilon_k}{2}\right)^2|f|^2(r)\right]dr\\
 &\ge -\frac{\epsilon_k}{2} \int_0^\infty \left(f^2\right)'(r)dr
 +\left(\frac{\epsilon_k}{2})\right)^2\int_0^\infty|f|^2(r)dr\\
 &\ge\frac{\epsilon_k}{2}|f|^2(0)+ \left(\frac{\epsilon_k}{2}\right)^2
 \int_0^\infty|f|^2(r)dr\\
 &\ge \left(\frac{\epsilon_k}{2}\right)^2\int_0^\infty |w|^2(r) e^{(n-1-2(k-1))r} dr\\
 &=\left(\frac{\epsilon_k}{2}\right)^2 \int_0^\infty |Mv|^2(r) e^{(n-1)r} dr.\\
 \end{split}
 \end{equation*}
 \hfill$\square$

 \noindent From the lemma (\ref{continuite}) and the estimate (\ref{decay2}) we obtain :
 \begin{equation}\label{decay3}
\|\beta_t\|_{L^2}\le \frac{2}{\frac{n-1}{2}-(k-1)}  \|\alpha\|_{L^2}.
 \end{equation}
 Hence if we let $t\to\infty$ in the equation (\ref{Cartan}), the estimates (\ref{decay1},\ref{decay3}) imply that
 that $\alpha=d\beta_\infty$ where 
 $$\beta_\infty=\int_0^\infty \left(\Phi^s\right)^*\left(\inte_T\alpha\right) ds\in L^2.$$
\endproof
\begin{rem}\label{primitive2rel}
When $k\ge (n+1)/2$, we can also show with a similar proof that
$\cH^k_{rel}(\Omega)=0$, even one can obtain that 
every $\alpha\in Z^k(\Omega,\partial \Omega)$ has an $L^2$ primitive given by
$$\beta=\int_{-\infty}^0(\Phi^s)^*\left(\inte_{T} \tilde\alpha\right)ds$$
where $\tilde \alpha$ is the extension of $\alpha$ to $\R\times N$ by letting 
$\tilde\alpha=0$ on $(\R\times N)\setminus \Omega=]-\infty,0[\times N$.
\end{rem}
\subsection{Mazzeo's result.}

We will now prove R.Mazzeo's result :
\begin{thm} Let $(M,g)$ be a conformally compact Riemannian manifold then
\begin{enumeroman}
\item when $k<\dim M/2$, then $H^k_2(M)\simeq H^k_0(M)\simeq H^k(\overline{M},\partial\overline{M} )$.
\item when $k>\dim M/2$, then $H^k_2(M)\simeq H^k(M)$.
\item when $k=\dim M/2$, then $\dim H^k_2(M)=\infty.$
\end{enumeroman}
\end{thm}
\proof We have already prove iii), we will only prove i) and indicate how we can prove
ii) with similar arguments. We can assume that the metric is exact. So let $\Omega=]0,\infty[\times N$
being a neighborhood of infinity endowed with the warped product metric
$(dr)^2+e^{2r} h$. And we will denote $K=M\setminus \Omega$ and $j_\Omega\,:\,\Omega\rightarrow M$,
 $j_K\,:\,K\rightarrow M$ the inclusion maps.

We assume that $k<n/2$. We are going to prove that the natural map
$H^k_0(M)\simeq H^k(K,\partial K)\rightarrow H^k_2(M)$ is an isomorphism.

{\bf fact 1:} The map $H^k_0(M)\simeq H^k(K,\partial K)\rightarrow H^k_2(M)$ is injective.
Let $\alpha\in C^\infty_0(\Lambda^k T^*M)$ a closed form with  support in $K$ which is mapped
 to zero in $H^k_2(M)$. According to
 the corollary (\ref{primitive})
 we know that $\alpha$ has a $L^2$ primitive : there is $\beta \in L^2(\Lambda^{k-1} T^*M)$
 such that 
 $$\alpha=d\beta .$$ Let $j_\Omega \,:\, \Omega\rightarrow M$ be the inclusion map, then
 we have that 
$$j_\Omega^*\beta\in Z^{k-1}_2(\Omega)$$ and by the proposition (\ref{primitive2}), we obtain 
$\eta\in L^2(\Lambda^{k-2} T^*\Omega)$ such that
$j_\Omega^*\beta=d\eta$. Consider $\bar \eta$ an extension of $\eta$ to $M$.
 Then $\beta-d\bar \eta$ has compact support and
$d(\beta-d\bar \eta)=\alpha$.

{\bf fact 2:} The map $H^k_0(M)\simeq H^k(K,\partial K)\rightarrow H^k_2(M)$ is surjective.
Let $\alpha\in \cH^k(M)$ then we know that $j^*_\Omega(\alpha)\in Z^k_2(\Omega)$ hence by 
the proposition (\ref{primitive2}), we obtain a
$\eta\in L^2(\Lambda^{k-1} T^*\Omega)$ such that
$d\eta=\alpha$.
Now we choose $\bar \eta\in L^2((\Lambda^{k-1} T^*M)$ an extension of
$\eta$ which is in the domain of $d$ and we obtain that
$\alpha-d\bar\eta$ has compact support moreover
 $\alpha$ and $\alpha-d\bar\eta$ belong to the same reduced cohomology class

The proof of the case ii) is done similarly. So assume that $k>\dim M/2$. 
First we recall that according to (\ref{vanishrel}), we have  $H^k_2(\Omega,\partial\Omega)=\{0\}$
hence the natural map 
$[j^*_K]\,:\, H^k_2(M)\rightarrow H^k(K)\simeq H^k(M)$ is injective (see \ref{exactrel}). 

In order to show that this map is surjective, we proceed as follow:
let $c\in H^k(K)$ and $\alpha\in c$  if $\bar \alpha\in C^\infty_0(\Lambda^kT^*M)$ is a smooth extension
of $\alpha$ then clearly
$$j_\Omega^*(d\bar\alpha)\in Z^{k+1}(\Omega,\partial\Omega).$$
We are then able to find some $\eta\in L^2(\Lambda^kT^*\Omega)$ which is smooth such that
$j_\Omega^*(d\bar\alpha)=d\eta$ and $\iota^*\eta=0$ 
(where $\iota\,:\, \partial\Omega\rightarrow \Omega$ is the inclusion map).
Then as in the proof of the lemma \ref{exactabs}, we can show that 
the $k$ form $\tilde \alpha$ defined
by $\tilde\alpha=\alpha$ on $K$ and 
$\tilde\alpha=\bar\alpha-\eta$ is closed and $L^2$ on $M$. Moreover it is clear that its $L^2$ cohomology
class
is map to $c$ by $[j^*_K]$.

When $k>(n-1)/2$, the surjectivity of this map is more easy to prove :
consider $\pi\,:\, M\rightarrow K$ the natural retraction which is the identity on $K$ and which is
$\pi(r,\theta)=(0,\theta)$ on $\Omega$. Then a small verification shows that
$\pi^*\alpha$ is a $L^2$ closed form on $M$ and clearly $j^*_K(\pi^*\alpha)=\alpha.$
\endproof
\begin{rem}
In fact, similar arguments shows that when $(M,g)$ is a complete Riemannian manifold such that 
$0$ is not in the spectrum of the Hodge-deRham Laplacian on $k$ forms then for $K\subset M$ a compact subset
and $\Omega=M\subset K$, we have the two short exact sequences
:
$$H_2^{k-1}(\Omega)\rightarrow H^k(K,\partial K)\rightarrow H_2^k(M)\rightarrow H_2^k(\Omega)\rightarrow H^{k+1}(K,\partial
K)\rightarrow H_2^{k+1}(M)\rightarrow H_2^{k+1}(\Omega). $$

$$H_2^{k-1}(K)\rightarrow H^k(\Omega,\partial \Omega)\rightarrow H_2^k(M)\rightarrow H_2^k(K)\rightarrow 
H^{k+1}(\Omega,\partial \Omega)\rightarrow H_2^{k+1}(M)\rightarrow H_2^{k+1}(K). $$

In fact, these exact sequences always hold for the non reduced $L^2$ cohomology and 
when $0$ is not in the spectrum of the Hodge-deRham Laplacian on $k$ forms then 
the non reduced and reduced $L^2$ cohomology coincide in degre $k$ and $k-1$ (this is exactly the assertion of the proposition
(\ref{primitive}).
\end{rem}
\subsection{Bibliographical hints}

The last lectures is based on a new proof of R. Mazzeo's results \cite{M} by N.
Yeganefar \cite{Y1}, some part of our arguments can be found in the paper
 of V.M. Gol'dshte\u\i n, V. I. Kuz'minov, I. A. Shvedov \cite{GKS}. A very good discussion on the essential spectrum can be find 
in the papers of H. Donnely \cite{Do} and C. B\"ar \cite{Ba} .  The integration by part
formula is due to H.Donnelly and F.Xavier \cite{DX} and has been revisited by many authors
for instance by J.Escobar-A.Freire \cite{EF} and A. Kasue \cite{Kasue} .

\end{document}